\documentclass[12pt]{article}
\usepackage{latexsym,amsmath,amssymb,amscd}
\begin{document}

\title{On Spherically Symmetric Motions of the Atmosphere Surrounding a Planet
Governed by the Compressible Euler Equations}
\author{Tetu Makino}
\date{\today}
\maketitle

\begin{abstract}
We consider spherically symmetric motions of inviscid compressible gas surrounding
a solid ball under the gravity of the core. Equilibria touch the vacuum with finite radii, and 
the linearized equation around one of the equilibria 
has time-periodic solutions. To justify the linearization, we should construct true solutions
for which this time-periodic solution plus the equilibrium is the first approximation. But this 
leads us to difficulty caused by singularities at the free boundary touching the
vacuum. We solve this problem
by the Nash-Moser theorem.


\end{abstract}

\newtheorem{Lemma}{Lemma}
\newtheorem{Proposition}{Proposition}
\newtheorem{Theorem}{Theorem}
\newtheorem{df}{Definition}
\section{Introduction}

We consider spherically symmetric motions of atmosphere governed by the 
compressible Euler equations:
\begin{align}
&\frac{\partial\rho}{\partial t} +u\frac{\partial \rho}{\partial r}+\rho\frac{\partial u}{\partial r}
+\frac{2}{r}\rho u =0, \nonumber \\
&\rho\Big(\frac{\partial u}{\partial t}+
u\frac{\partial u}{\partial r}\Big)+\frac{\partial P}{\partial r}=-\frac{g_0\rho}{r^2}\quad
(R_0\leq r )
\end{align}
and the boundary value condition
\begin{equation}
\rho u|_{r=R_0}=0.
\end{equation}
Here $\rho$ is the density, $u$ the velocity, $P$ the pressure.
$R_0\ (>0)$ is the radius of the central solid ball, and $g_0=G_0M_0$, $G_0$ being the gravitational constant,
$M_0$ the mass of the central ball. The self-gravity of the atmosphere is neglected.

In this study we always assume that 
\begin{equation}
P=A\rho^{\gamma},
\end{equation}
where $A$ and $\gamma$ are positive constants, and we assume that $1<\gamma \leq 2$. \\

Equilibria of the problem are given by
$$\bar{\rho}(r)=\begin{cases}
\displaystyle A_1\Big(\frac{1}{r}-\frac{1}{R}\Big)^{\frac{1}{\gamma-1}} & \quad (R_0\leq r<R) \\
0 & \quad (R\leq r),
\end{cases}
$$
where $R$ is an arbitrary number such that $R>R_0$ and
$$A_1=\Big(\frac{(\gamma-1)g_0}{\gamma A}\Big)^{\frac{1}{\gamma-1}}.$$

{\bf Remark} The total mass $M$ of the equilibrium is given by
$$M=4\pi A_1\int_{R_0}^R\Big(\frac{1}{r}-\frac{1}{R}\Big)^{\frac{1}{\gamma-1}}r^2dr.$$
$M$ is an increasing function of $R$. Of course
$M\rightarrow 0$ as $R \rightarrow R_0$. But as $R\rightarrow +\infty$, we see
$$M\rightarrow 4\pi A_1\int_{R_0}^{\infty}r^{\frac{2\gamma-3}{\gamma-1}}dr
=\begin{cases}
+\infty & \mbox{if}\  \gamma\geq 4/3 \\
M^*(<\infty) & \mbox{if}\  \gamma <4/3,
\end{cases}
$$
where
$$M^*=\frac{4\pi A_1(\gamma-1)}{4-3\gamma}R_0^{-\frac{4-3\gamma}{\gamma-1}}.$$
Hence if $\gamma \geq 4/3$ there is an equilibrium for
any given total mass, but
if $\gamma <4/3$ the possible mass has the upper bound
 $M^*$. 
Anyway, given the total mass $M$, a conserved quantity, in $(0,+\infty)$
or $(0,M^*)$, then the radius $R$ or the configuration of the equilibrium
is uniquely determined.\hfill$\square$

Let us fix one of these equilibria. We are interested in motions around
this equilibrium. \\

Here let us glance at the history of researches of this problem.

Of course there were a lot of works on the Cauchy problem to the
compressible Euler equations. But there were gaps if we consider 
density distributions which contain vacuum regions.

As for local-in-time existence of smooth density
with compact support, \cite{M1989} treated the problem under the assumption
that the initial density is non-negative and the initial value of
$$\omega:=\frac{2\sqrt{A\gamma}}{\gamma-1}\rho^{\frac{\gamma-1}{2}}$$
is smooth, too. By the variables $(\omega,u)$ the equations are symmetrizable
continuously including the region of vacuum. Hence the theory of 
quasi-linear symmetric hyperbolic
systems can be applied. The discovery of the variable $\omega$ can go back to
\cite{M1986}, \cite{MUK}.
However, since
$$\omega\propto \Big(\frac{1}{r}-\frac{1}{R}\Big)^{\frac{1}{2}}\sim 
\mbox{Const.}(R-r)^{\frac{1}{2}} \quad \mbox{as}\  r\rightarrow R-0$$
for equilibria, $\omega$ is not smooth at the boundary $r=R$ with the vacuum.
Hence the class of ``tame" solutions considered in \cite{M1989} cannot cover equilibria.

On the other hand, possibly discontinuous weak solutions with compactly
supported density can be constructed. The article \cite{MT} gave local-in-time
existence of bounded weak solutions under the assumption
that the initial density is bounded and non-negative. The proof by the 
compensated compactness method is due to
\cite{DCL}. Of course the class of weak solutions can cover equilibria, 
but the concrete structures of solutions were not so clear. 

Therefore we wish to construct solutions whose regularities 
are weaker than solutions with smooth $\omega$ and stronger than possibly discontinuous 
weak solutions. The present result is an answer to this wish. 
More concretely speaking, the solution $(\rho(t,r),u(t,r))$ constructed
in this article should be continuous on
$0\leq t\leq T,R_0\leq r <\infty$ and there should be found a continuous curve
$r=R_F(t), 0\leq t\leq T,$ such that
$|R_F(t)-R|\ll 1, \rho(t,r)>0 $ for $  0\leq t\leq T, R_0\leq r <R_F(t)$
and $\rho(t,r)=0$ for $0\leq t\leq T, R_F(t)\leq r<\infty$. The curve $r=R_F(t)$ is
the free boundary at which the density touches the vacuum.
It will be shown that the solution satisfies
$$\rho(t,r)=C(t)(R_F(t)-r)^{\frac{1}{\gamma-1}}(1+O(R_F(t)-r))$$
as $r \rightarrow R_F(t)-0$. Here $C(t)$ is positive and smooth in $t$.
This situation is ``physical vacuum boundary" so-called by
\cite{JM} and \cite{CS}. This concept can be traced back to
\cite{L}, \cite{LY}, \cite{Y}. Of course this singularity is just that
of equilibria.\\

The major difficulty of the analysis
comes from the free boundary touching the vacuum, which can move along time. 
So it is convenient to
introduce the Lagrangian
mass coordinate
$$m=4\pi\int_{R_0}^r\rho(t,r')r'^2dr',$$
to fix the interval of independent variable to consider.
Taking $m$ as the independent variable instead of $r$,
the equations turn out to be
\begin{align*}
&\frac{\partial\rho}{\partial t}+4\pi \rho^2
\frac{\partial}{\partial m}(r^2u)=0, \\
&\frac{\partial u}{\partial t}+
4\pi r^2\frac{\partial P}{\partial m}=-\frac{g_0}{r^2} \qquad (0<m<M),
\end{align*}
where
$$r=\Big(R_0^3+\frac{3}{4\pi}\int_0^m\frac{dm}{\rho}\Big)^{1/3}.$$

We note that
$$\frac{\partial r}{\partial t}=u, \qquad \frac{\partial r}{\partial m}=\frac{1}{4\pi r^2\rho}.$$

Let us take $\bar{r}=\bar{r}(m)$ as the independent variable instead of $m$, where
$m\mapsto \bar{r}=\bar{r}(m)$ is the inverse function of the function
$$\bar{r}\mapsto m=4\pi \int_{R_0}^{\bar{r}}\bar{\rho}(r')r'^2dr'.$$
Then, since 
$$\frac{\partial }{\partial m}=\frac{1}{4\pi\bar{r}^2\bar{\rho}}\frac{\partial}{\partial\bar{r}},
\qquad
\rho=\Big(4\pi r^2\frac{\partial r}{\partial m}\Big)^{-1}=\bar{\rho}
\Big(\frac{r^2}{\bar{r}^2}\frac{\partial r}{\partial\bar{r}}\Big)^{-1},
$$
we have a single second-order equation
$$
\frac{\partial^2r}{\partial t^2}+
\frac{1}{\bar{\rho}}
\frac{r^2}{\bar{r}^2}\frac{\partial}{\partial\bar{r}}
\Big(\bar{P}\Big(\frac{r^2}{\bar{r}^2}\frac{\partial r}{\partial\bar{r}}\Big)^{-\gamma}\Big)+
\frac{g_0}{r^2}=0.
$$
The variable $\bar{r}$ runs on the interval $[R_0, R]$ and
the boundary condition is
$r|_{\bar{r}=R_0}=R_0.$ 

\textbullet {\it Without loss of generality, we can and shall assume that}
$$R_0=1,\quad g_0=\frac{1}{\gamma-1}, \quad A=\frac{1}{\gamma}, \quad A_1=1.$$

Keeping in mind 
that the equilibrium satisfies
$$\frac{1}{\bar{\rho}}\frac{\partial \bar{P}}{\partial\bar{r}}+\frac{g_0}{\bar{r}^2}=0, $$
we have
$$
\frac{\partial^2r}{\partial t^2}-
\frac{1}{\bar{\rho}}
\frac{r^2}{\bar{r}^2}\frac{\partial}{\partial\bar{r}}
\Big(\bar{P}\Big(1-
\Big(\frac{r^2}{\bar{r}^2}\frac{\partial r}{\partial\bar{r}}\Big)^{-\gamma}\Big)\Big)+
\frac{1}{\gamma-1}\Big(\frac{1}{r^2}-\frac{r^2}{\bar{r}^4}\Big)=0.$$
Introducing the unknown variable $y$ for perturbation by
\begin{equation}
r=\bar{r}(1+y),
\end{equation}
we can write the equation as
\begin{equation}
\frac{\partial^2y}{\partial t^2}-
\frac{1}{\rho r}
(1+y)^2\frac{\partial}{\partial r}\Big(PG\Big(y, r\frac{\partial y}{\partial r}\Big)\Big)
-\frac{1}{\gamma-1}\frac{1}{r^3}H(y)=0,
\end{equation}
where
\begin{align*}
G(y,v)&:=1-(1+y)^{-2\gamma}(1+y+v)^{-\gamma}=\gamma(3y+v)+[y,v]_2, \\
H(y)&:=(1+y)^2-\frac{1}{(1+y)^2}=4y+[y]_2
\end{align*}
and we have used the abbreviations $r, \rho, P$
for $\bar{r}, \bar{\rho}, \bar{P}$.

{\bf Notational Remark} Here and hereafter $[X]_q$ denotes a convergent power series, 
or an analytic function given by the series, of the form $\sum_{j\geq q}a_jX^j$, and
$[X,Y]_q$ stands for a convergent double power series of the form
$\sum_{j+k\geq q}a_{jk}X^jY^k$.\hfill$\square$

We are going to study the equation (5) on $1<r<R$ with the boundary condition
$$y|_{r=1}=0.$$
Of course $y$ and $\displaystyle r\frac{\partial y}{\partial r}$ will be confined to 
$$|y|+\Big|r\frac{\partial y}{\partial r}\Big| <1.$$

Here let us propose the main goal of this study roughly.
Let us fix an arbitrarily large positive number $T$. Then we have \\

{\bf Main Goal } {\it For sufficiently small $\varepsilon>0$ there
is a solution $y=y(t,r;\varepsilon)$ of (5) in
$C^{\infty}([0,T]\times[1,R])$ such that
$$y(t,r;\varepsilon)=\varepsilon y_1(t,r)+O(\varepsilon^2).$$
The same estimates $O(\varepsilon^2)$ hold between the
higher order derivatives of $y$ and $\varepsilon y_1$.}\\

Here $y_1(t,r)$ is a time-periodic function specified in Section 2, which
is of the form
$$y_1(t,r)=\sin(\sqrt{\lambda}t+\theta_0)\cdot \tilde{\Phi}(r),$$
where $\lambda$ is a positive number, $\theta_0$ a constant, and
$\tilde{\Phi}(r)$ is an analytic function of $1\leq r\leq R$.

Once the solution $y(t,r;\varepsilon)$ is given, then
 the corresponding motion of gas particles can be expressed by the Lagrangian 
coordinate as
\begin{align*}
r(t,m)&=\bar{r}(m)(1+y(t,\bar{r}(m);\varepsilon)) \\
&=\bar{r}(m)(1+\varepsilon y_1(t,\bar{r}(m))+O(\varepsilon^2)).
\end{align*}
The curve $r=R_F(t)$ of the free vacuum boundary is given by
$$R_F(t)=r(t,M)=R+\varepsilon R\sin(\sqrt{\lambda}t+\theta_0)\tilde{\Phi}(R)+O(\varepsilon^2).$$
{\it The free boundary $R_F(t)$ oscillates around $R$ with time-period $2\pi/\sqrt{\lambda}$ approximately.}

The solution $(\rho,u)$ of the original problem (1)(2)
is given by
$$\rho=\bar{\rho}(\bar{r})\Big((1+y)^2\Big(1+
y+\bar{r}\frac{\partial y}{\partial\bar{r}}\Big)\Big)^{-1},
\qquad u=\bar{r}\frac{\partial y}{\partial t} $$
implicitly by
\begin{align*}
\bar{r}&=\bar{r}(m), \qquad y=y(t,\bar{r}(m);\varepsilon) \\
\frac{\partial y}{\partial\bar{r}}&=
\partial_ry(t,\bar{r}(m);\varepsilon), \qquad
\frac{\partial y}{\partial t}=
\partial_ty(t,\bar{r}(m);\varepsilon),
\end{align*}
where $m=m(t,r)$ for $1\leq r\leq R_F(t)$. Here
$r \mapsto m=m(t,r)$ is the inverse function of the 
function
$m\mapsto r=r(t,m)=\bar{r}(m)(1+y(t,\bar{r}(m);\varepsilon))$.
We note that
$$R_F(t)-r(t,m)=R(1+y(t,R))-\bar{r}(m)(1+y(t,\bar{r}(m))$$
implies
$$\frac{1}{\kappa}(R-\bar{r})\leq R_F(t)-r\leq \kappa (R-\bar{r})
$$
with $0<\kappa-1\ll 1$, since
$|y|+|\partial_r y|\leq \varepsilon C$. Therefore
$$y(t,\bar{r}(m))=y(t,R)+O(R_F(t)-r),$$
and so on.
Hence we get the ``physical vacuum boundary". In fact
the corresponding density distribution 
$\rho=\rho(t,r)$, where $r$ is the original Euler coordinate, satisfies
$$\rho(t,r)>0\  \mbox{for}\  1\leq r<R_F(t), \qquad
\rho(t,r)=0\  \mbox{for}\  R_F(t)\leq r, $$
and, since $y(t,r)$ is smooth on $1\leq r\leq R$, we have
$$\rho(t,r)=C(t)(R_F(t)-r)^{\frac{1}{\gamma-1}}(1+O(R_F(t)-r)) $$
as $r\rightarrow R_F(t)-0$. Here $C(t)$ is positive and smooth in $t$.\\

We shall give a precise statement of the main result in Section 3
and give a proof of the main result in Sections 4, 5.
We shall apply
the Nash-Moser theory. The reason is as follows.

The equation (5) looks like as if it is a second-order quasi-linear
hyperbolic equation, and one might expect that the usual
iteration method in a suitable Sobolev spaces, e.g., $H^s$, or 
something like them, could be used. But it is not the case. Actually the linear part of the 
equation is essentially
the d'Alembertian operator
$$
\frac{\partial^2}{\partial t^2}-\triangle=
\frac{\partial^2}{\partial t^2}-x\frac{\partial^2}{\partial x^2}-\frac{N}{2}\frac{\partial}{\partial x}
=\frac{\partial^2}{\partial t^2}- \frac{\partial^2}{\partial\xi^2}-\frac{N-1}{\xi}\frac{\partial}{\partial\xi}
$$
in the variables $x, \xi$ such that 
$\displaystyle \frac{R-r}{R} \sim x=\frac{\xi^2}{4}, $
and
the nonlinear terms are smooth functions of $y$ and $\partial y/\partial r$. (See (13) and (15).)
Here the term $\partial y/\partial r$ apparently looks like as if to be the first order derivative.
If it was the case, the usual iteration in the Picard's scheme applied to the wave
equation would work, since the inverse of the d'Alembertian operator
may recover the regularity up to one order of derivative,
that is, roughly speaking, the d'Alembertian may
pull back $C^1([0,T],L^2)$ to $C^1([0,T],H^1)$. But indeed
the apparently first order derivative
$\partial y/\partial r$ performs like
$$r\frac{\partial y}{\partial r} 
\sim \frac{\partial y}{\partial x}\propto -\frac{1}{\xi}\frac{\partial y}{\partial \xi}
\sim -\frac{\partial^2y}{\partial\xi^2}$$
near $r=R$ or $\xi=0$, like the second order derivative.
So, since the inverse of the d'Alembertian recovers the regularity up to only
one order of derivative, the usual iteration for quasi-linear wave equations
may cause troubles of the loss of regularities, which occur at the free
vacuum boundary $r=R$. This is the reason why we like to apply 
the Nash-Moser theory to our problem. 

In fact, roughly speaking, the essence of
the Nash-Moser method is as follows:
Suppose that we want to solve the equation $\mathfrak{P}(w)=y$, where
$\mathfrak{P}$ is a nonlinear differential operator of the second order, that is,
if $w$ is $n+2$ times differentiable, then $y$ is $n$ times differentiable, assuming $\mathfrak{P}(0)=0$
; Suppose 
the linearized equation $D\mathfrak{P}(w)h=g$ admits $n+2-l$ times differentiable
solutions $h$ for given $n$ times differentiable
$g$ when $n+2$ times differentiable $w$ is fixed;
Here $$D\mathfrak{P}(w)h:=\lim_{\tau\rightarrow 0}
\frac{1}{\tau}(\mathfrak{P}(w+\tau h)-\mathfrak{P}(w));$$ If $l>0$, the simple 
Picard's iteration
$$w_{\nu+1}=w_{\nu}-D\mathfrak{P}(0)^{-1}(\mathfrak{P}(w_{\nu})-y)$$
does not work, but, even if
$l>0$, that is, even if there is a `derivative loss', the equation can be solved
by a Newton's approximation combined with so called
`smoothing operators'. (If $D\mathfrak{P}(w)$ is essentially
d'Alembertian, we can consider $l=1$. )
Various formulations of the Nash-Moser theory are known.
Among them the present study will adopt the formulation by R. Hamilton 
in \cite{Hamilton}. See Section 4.

\section{Analysis of the linearized problem}

The linearization of the equation (5) is clearly
\begin{equation}
\frac{\partial^2y}{\partial t^2}+\mathcal{L}\Big(r,\frac{\partial}{\partial r}\Big)y=0,
\end{equation}
where
\begin{align}
&\mathcal{L}y=\mathcal{L}\Big(r,\frac{d}{dr}\Big)y :=-\frac{1}{\rho r}\frac{d}{dr}
\Big(P\gamma\Big(3y+r\frac{dy}{dr}\Big)\Big)-\frac{1}{\gamma-1}\frac{1}{r^3}(4y) \nonumber \\
&=-\Big(\frac{1}{r}-\frac{1}{R}\Big)\frac{d^2y}{dr^2}+
\Big(-\frac{4}{r}\Big(\frac{1}{r}-\frac{1}{R}\Big)+
\frac{\gamma}{\gamma-1}\frac{1}{r^2}\Big)\frac{dy}{dr}+\frac{3\gamma-4}{\gamma-1}\frac{y}{r^3}.
\end{align}

In order to analyze the eigenvalue problem
$\mathcal{L} 
y=\lambda y$, we
introduce the independent variable $z$ by
\begin{equation}
z=\frac{R-r}{R}
\end{equation}
and the parameter $N$ by
\begin{equation}
\frac{\gamma}{\gamma-1}=\frac{N}{2} \quad \mbox{or} \quad \gamma=1+\frac{2}{N-2}.
\end{equation}

Then we can write
\begin{equation}
R^3\mathcal{L}y=-\frac{z}{1-z}\frac{d^2y}{dz^2}
-\frac{\frac{N}{2}-4z}{(1-z)^2}\frac{dy}{dz}+\frac{8-N}{2}\frac{1}{(1-z)^3}y.
\end{equation}

The variable $z$ runs over the interval $[0, 1-1/R]$, the boundary $z=0$ corresponds to the free
boundary touching the vacuum, and the boundary condition at $z=1-1/R$
is the Dirichlet condition $y=0$. 

Although the boundary $z=1-1/R$ is regular, the boundary $z=0$ is singular.
In order to analyze the singularity,
we transform the equation $\mathcal{L}y=\lambda y$ , which can be written as
$$
-z\frac{d^2y}{dz^2}-\Big(\frac{N}{2}\frac{1}{1-z}-\frac{4z}{1-z}\Big)\frac{dy}{dz}
+\frac{8-N}{2}\frac{y}{(1-z)^2}=\lambda R^3(1-z)y,
$$
to an equation of the formally self-adjoint
form
$$-\frac{d}{dz}p(z)\frac{dy}{dz}+q(z)y=\lambda R^3\mu(z)y.$$
This can be done by putting
$$
p=z^{\frac{N}{2}}(1-z)^{\frac{8-N}{2}}, \quad
q=\frac{8-N}{2}z^{\frac{N-2}{2}}(1-z)^{\frac{4-N}{2}}, \quad
\mu=z^{\frac{N-2}{2}}(1-z)^{\frac{10-N}{2}}.
$$
Using the Liouville transformation, we convert the equation
$$-\frac{d}{dz}p(z)\frac{dy}{dz}+q(z)y=\lambda R^3\mu(z)y+f$$
to the standard form
$$-\frac{d^2\eta}{d\xi^2}+Q\eta =\lambda R^3\eta +\hat{f}.$$
This can be done by putting
\begin{align}
\xi&=\int_0^z\sqrt{\frac{\mu}{p}}dz=\int_0^z\sqrt{\frac{1-\zeta}{\zeta}}d\zeta 
=\sqrt{z(1-z)}+\tan^{-1}\sqrt{\frac{z}{1-z}}, \\
\eta&=(\mu p)^{1/4}y=z^{\frac{N-1}{4}}(1-z)^{\frac{9-N}{4}}y,\quad
\hat{f}=p^{1/4}\mu^{-3/4}f=z^{\frac{3-N}{4}}(1-z)^{\frac{N-11}{4}}f, \nonumber
\end{align}
and
\begin{align*}
Q&=\frac{p}{\mu}\Big(\frac{q}{p}+
\frac{1}{4}\Big(\frac{p'}{p}+\frac{\mu'}{\mu}\Big)'
-\frac{1}{16}\Big(\frac{p'}{p}+\frac{\mu'}{\mu}\Big)^2
+\frac{1}{4}\frac{p'}{p}\Big(\frac{p'}{p}+\frac{\mu'}{\mu}\Big)\Big) \\
&=\frac{1}{z(1-z)^3}\Big(\frac{(N-1)(N-3)}{16}+\frac{7-2N}{2}z+2z^2\Big).
\end{align*}
Putting 
$$\xi_R :=\int_0^{1-\frac{1}{R}}\sqrt{\frac{1-z}{z}}dz \quad (<\frac{\pi}{2}), $$
we see that the variable $\xi$ runs over the interval $[0, \xi_R]$.
Since $z\sim \displaystyle \frac{\xi^2}{4}$ as $ \xi \rightarrow 0$, we see
$$Q \sim \frac{(N-1)(N-3)}{4}\frac{1}{\xi^2}$$
as $\xi \rightarrow 0$. But $\gamma <2$ implies $N>4$ and
$\displaystyle \frac{(N-1)(N-3)}{4}>\frac{3}{4}$. Hence the boundary $\xi=0$ is of the limit point type.
See, e.g., \cite[p.159, Theorem X.10.]{Reed} 
The exceptional case $\gamma=2$ or $N=4$ will be considered separately.
Anyway the potential $Q$ is bounded from below on
$0<\xi<\xi_R$ provided that $N>3$. 
Thus we have

\begin{Proposition}
The operator $T_0, \mathcal{D}(T_0)=C_0^{\infty}(0,\xi_R), T_0\eta=
-\eta_{\xi\xi}+Q\eta$, in $L^2(0,\xi_R)$ has the Friedrichs extension
$T$, a self-adjoint operator whose spectrum consists of simple
eigenvalues
$\lambda_1R^3<\lambda_2R^3<\cdots<\lambda_nR^3<\cdots\rightarrow +\infty.$
In other words, the operator
$S_0, \mathcal{D}(S_0)=C_0^{\infty}(0, 1-\frac{1}{R}), S_0y=\mathcal{L}y$
in 
$$\mathfrak{X}:=L^2((0,1-\frac{1}{R}), \mu dz(=z^{\frac{N-2}{2}}(1-z)^{\frac{10-N}{2}}dz)),$$
 has
the Friedrichs extension $S$, a self-adjoint operator with the eigenvalues $(\lambda_n)_n$ .
\end{Proposition}

We note that the domain of $S$ is
\begin{align*}
\mathcal{D}(S)&=\{ y \in \mathfrak{X} \ | \ \exists \phi_n \in C_0^{\infty}(0, 1-\frac{1}{R}) \\
&\mbox{such that}\  \phi_n\rightarrow y \ \mbox{in} \ \mathfrak{X} \ \mbox{and}\  
\mathfrak{Q}[\phi_m-\phi_n]\rightarrow 0 \ \mbox{as}\ m,n\rightarrow \infty, \\
&\mbox{and}\  \mathcal{L}y\in \mathfrak{X} \ \mbox{in distribution sense}. \}
\end{align*}
Here
$$
\mathfrak{Q}[\phi]:=\int_0^{1-\frac{1}{R}}\Big|\frac{d\phi}{dz}\Big|^2z^{\frac{N}{2}}(1-z)^{\frac{8-N}{2}}dz 
=\int_0^{1-\frac{1}{R}}\frac{z}{1-z}\Big|\frac{d\phi}{dz}\Big|^2\mu(z)dz.
$$

Moreover we have

\begin{Proposition}
If $N\leq 8$ ( or $\gamma\geq 4/3$ ), the least eigenvalue $\lambda_1$ is positive.
\end{Proposition}
{\bf Proof}
Suppose $N\leq 8$. Clearly $y \equiv 1$ satisfies
$$-\frac{d}{dz}p\frac{dy}{dz}+qy=q=\frac{8-N}{2}z^{\frac{N-2}{2}}(1-z)^{\frac{4-N}{2}}.$$
Therefore the corresponding $\eta_1(\xi)$ given by
$$\eta_1=z^{\frac{N-1}{4}}(1-z)^{\frac{9-N}{4}}$$
satisfies
$$-\frac{d^2\eta_1}{d\xi^2}+Q\eta_1=\hat{q}=\frac{8-N}{2}z^{\frac{N-1}{4}}(1-z)^{-\frac{N+3}{4}}.$$
It is easy to see $\eta_1=\displaystyle \frac{d\eta_1}{d\xi}=0$ at $\xi=0$. Let
$\phi_1(\xi)$ be the eigenfunction of
$-d^2/d\xi^2 +Q$ associated with the least eigenvalue $\lambda_1$.
We can assume that $\phi_1(\xi)>0$ for $0<\xi<\xi_R, \phi_1(\xi_R)=0$, 
and
$\displaystyle \frac{d\phi_1}{d\xi}<0$ at $\xi=\xi_R$.
Then integrations by parts give
\begin{align*}
&\lambda_1\int_0^{\xi_R}\phi_1\eta_1d\xi =
\int_0^{\xi_R}\Big(-\frac{d^2\phi_1}{d\xi^2}+Q\phi_1\Big)\eta_1d\xi \\
&=-\frac{d\phi_1}{d\xi}\eta_1\Big|_{\xi=\xi_R}+
\int_0^{\xi_R}\Big(\frac{d\phi_1}{d\xi}\frac{d\eta_1}{d\xi}+Q\phi_1\eta_1\Big)d\xi 
>\int_0^{\xi_R}\Big(\frac{d\phi_1}{d\xi}\frac{d\eta_1}{d\xi}+Q\phi_1\eta_1\Big)d\xi \\
&=\int_0^{\xi_R}\phi_1\Big(-\frac{d^2\eta_1}{d\xi^2}+Q\eta_1\Big)d\xi 
=\int_0^{\xi_R}\phi_1\hat{q}d\xi \geq 0.
\end{align*}
\hfill$\blacksquare$

{\bf Remark} When $N=8$, $\eta_1$ satisfies
$\displaystyle -\frac{d^2\eta}{d\xi^2}+Q\eta=0$, but does not satisfy the
boundary condition $\eta|_{\xi=\xi_R}=0$. Hence
it is not an eigenfunction of zero eigenvalue, and $\lambda_1>0$ even if $N=8$.\hfill$\square$

For the sake of convenience of the further analysis,
let us rewrite the linear part $\mathcal{L}$ by introducing a new variable
$$\tilde{x}=\frac{\xi^2}{4}=\frac{1}{4}\Big(\sqrt{z(1-z)}+\tan^{-1}\sqrt{\frac{z}{1-z}}\Big)^2.
$$
Clearly $\tilde{x}=z+[z]_2$ and the change of variables
$z \mapsto \tilde{x}$ is analytic on $0\leq z <1$ and its inverse
$\tilde{x} \mapsto z$ is analytic on $0\leq \tilde{x} <\tilde{x}_{\infty}:=\pi^2/16$.
Since
\begin{align*}
\frac{d}{dz}&=\sqrt{\frac{\tilde{x}}{z}}\sqrt{1-z}\frac{d}{d\tilde{x}}, 
\quad\mbox{with}\quad\sqrt{\frac{\tilde{x}}{z}}=1+[\tilde{x}]_1, \\
\frac{d^2}{dz^2}&=\frac{1-z}{z}\Big(
\tilde{x}\frac{d^2}{d\tilde{x}^2}+
\frac{1}{2}\Big(1-\sqrt{\frac{\tilde{x}}{z}}\frac{1}{(1-z)\sqrt{1-z}}\Big)\frac{d}{d\tilde{x}}\Big),
\end{align*}
we can write
$$R^3\mathcal{L}y=
-\Big(\tilde{x}\frac{d^2y}{d\tilde{x}^2}+\frac{N}{2}\frac{dy}{d\tilde{x}}\Big)+
\ell_1(\tilde{x})\tilde{x}\frac{dy}{d\tilde{x}}+
\ell_0(\tilde{x})y,
$$
where
$\ell_1(\tilde{x})$ and $\ell_0(\tilde{x})$ are analytic on $0\leq \tilde{x}<\tilde{x}_{\infty}$.
Putting
\begin{align}
x&=R^3\tilde{x}=\frac{R^3\xi^2}{4}=\frac{R^3}{4}\Big(\sqrt{z(1-z)}+\tan^{-1}\sqrt{\frac{z}{1-z}}\Big)^2\nonumber \\
&=\frac{R^3}{4}\Big(\frac{\sqrt{(R-r)r}}{R}+\tan^{-1}\sqrt{\frac{R-r}{r}}\Big)^2,
\end{align}
we can write
\begin{equation}
\mathcal{L}y=-\triangle y+L_1(x)x\frac{dy}{dx}+L_0(x)y,
\end{equation}
where 
$$\triangle =x\frac{d^2}{dx^2}+\frac{N}{2}\frac{d}{dx}$$
and $L_1(x)$ and $L_0(x)$ are analytic on
$0\leq x < x_{\infty}:=R^3\tilde{x}_{\infty}=\pi^2R^3/16$. While $r$ runs over the interval
$[1,R]$, $x$ runs over $[0, x_R]$, where
$x_R:=R^3\xi_R^2/4 (< x_{\infty})$. The Dirichlet condition at the regular boundary is
$\displaystyle y|_{x=x_R}=0$.

{\bf Remark}
Since $x=R^3\xi^2/4$, we have
$$\triangle=x\frac{d^2}{dx^2}+\frac{N}{2}\frac{d}{dx}=
\frac{1}{R^3}\Big(\frac{d^2}{d\xi^2}+\frac{N-1}{\xi}\frac{d}{d\xi}\Big).$$
Thus $\triangle$ is the radial part of the Laplacian in the
$N$-dimensional Euclidean space $\mathbb{R}^N$ provided that $N$
is an integer. But we do not assume that
$N$ is an integer in this study.\hfill$\square$

Let us fix a positive eigenvalue $\lambda=\lambda_n$ and an associated
eigenfunction $\Phi(x)$ of $\mathcal{L}$. Then
\begin{equation}
y_1(t,x)=\sin(\sqrt{\lambda}t+\theta_0)\Phi(x)
\end{equation}
is a time-periodic solution of the linearized problem
$$\frac{\partial^2y}{\partial t^2}+\mathcal{L}y=0, \qquad y|_{x=x_R}=0. $$
Moreover we claim that $\Phi(x)$ is an analytic function of $0\leq x<x_{\infty}$. 
To verify it, we use the following

\begin{Lemma}
We consider the differential equation
$$x\frac{d^2y}{dx^2}+b(x)\frac{dy}{dx}+c(x)y=0,$$
where
$$b(x)=\beta+[x]_1, \qquad c(x)=[x]_0, $$
and we assume that $\beta \geq 2$. Then 1) there is a solution $y_1$
of the form
$$y_1=1+[x]_1, $$
and 2) there is a solution $y_2$ such that
$$y_2=x^{-\beta+1}(1+[x]_1)
$$
provided that $\beta \not\in \mathbb{N}$ or
$$y_2=x^{-\beta+1}(1+[x]_1)+hy_1\log x$$
provided that $\beta\in\mathbb{N}$. Here $h$ is a constant which can vanish in some cases.
\end{Lemma}

For a proof, see \cite[Chapter 4]{Coddington}. Applying this lemma with $\beta=N/2$ to the equation
$$
x\frac{d^2y}{dx^2}+\Big(\frac{N}{2}-L_1(x)x\Big)\frac{dy}{dx}
+(\lambda-L_0(x))y=0, $$
we get the assertion, since $y_2\sim \displaystyle
x^{-\frac{N-2}{2}}$ cannot belong to
$\mathfrak{X}=L^2(x^{\frac{N-2}{2}}dx)$ for $N\geq 4$, even if
$N=4$, which was the exceptional case in the preceding discussion of the limit point type.

\section{Statement of the main result}

%\subsection{}
We rewrite the equation (5) by using the linearized part
$\mathcal{L}$ defined by (7) as
\begin{equation}
\frac{\partial^2y}{\partial t^2}+
\Big(1+G_I\Big(y,r\frac{\partial y}{\partial r}\Big)\Big)\mathcal{L}y+
G_{II}\Big(r,y,r\frac{\partial y}{\partial r}\Big)=0,
\end{equation}
where
\begin{align*}
G_I(y,v)&=(1+y)^2\Big(1+\frac{1}{\gamma}\partial_vG_2(y,v)\Big)-1, \\
G_{II}(r,y,v)&=\frac{P}{\rho r^2}G_{II0}(y,v)+
\frac{1}{\gamma-1}\frac{1}{r^3}G_{II1}(y,v), \\
G_{II0}(y,v)&=(1+y)^2(3\partial_vG_2-\partial_yG_2)v, \\
G_{II1}(y,v)&=(1+y)^2\Big(
\frac{1}{\gamma}(\partial_vG_2)((4-3\gamma)y-\gamma v)+G_2\Big)
-H+4y(1+y)^2.
\end{align*}
Here
$$
G_2:=G-\gamma(3y+v)=[y,v]_2, $$
$$\partial_vG_2=\frac{\partial G}{\partial v}-\gamma=[y,v]_1,\quad
\partial_yG_2=\frac{\partial G}{\partial y}-3\gamma=[y,v]_1.
$$

We have fixed a solution $y_1$ of the linearized equation 
$y_{tt}+\mathcal{L}y=0$ (see (14)), and we seek a solution $y$ of (5) or (15) of the form
$$y=\varepsilon y_1+\varepsilon w,$$
where $\varepsilon$ is a small positive parameter. 

{\bf Remark} The following discussion is valid if we take
$$y_1=\sum_{k=1}^Kc_k\sin(\sqrt{\lambda_{n_k}}t+\theta_k)\cdot \Phi_k(x), \eqno(14)'$$
where $\Phi_k$ is an eigenfunction of $\mathcal{L}$ associated with the positive eigenvalue
$\lambda_{n_k}$ and $c_k$ and $\theta_k$ are constants for $k=1,\cdots, K$.\hfill$\square$

Then the equation which governs $w$
turns out to be
\begin{equation}
\frac{\partial^2w}{\partial t^2}+
\Big(1+\varepsilon a(t,r,w,r\frac{\partial w}{\partial r},\varepsilon)\Big)
\mathcal{L}w+
\varepsilon b(t,r,w,r\frac{\partial w}{\partial r},\varepsilon)=
\varepsilon c(t,r,\varepsilon),
\end{equation}
where
\begin{align*}
a(t,r,w, \Omega,\varepsilon)&=\varepsilon^{-1}G_I(\varepsilon(y_1+w),\varepsilon(v_1+\Omega)), \\
b(t,r,w,\Omega,\varepsilon)&=\varepsilon^{-1}G_I(\varepsilon(y_1+w),\varepsilon(v_1+\Omega))\mathcal{L}
y_1+\varepsilon^{-2}G_{II}
(r,\varepsilon(y_1+w),\varepsilon(v_1+\Omega)) \\
&-\varepsilon^{-1}G_I(\varepsilon y_1, \varepsilon v_1)\mathcal{L}y_1
-\varepsilon^{-2}G_{II}(r,\varepsilon y_1,\varepsilon v_1), \\
c(t,r,\varepsilon)&=-\varepsilon^{-1}G_I(\varepsilon y_1, \varepsilon v_1)\mathcal{L}y_1
-\varepsilon^{-2}G_{II}(r,\varepsilon y_1,\varepsilon v_1).
\end{align*}
Here $v_1$ stands for $r\partial y_1/\partial r$.

The main result of this study can be stated as follows:

\begin{Theorem}
For any $T>0$, there is a sufficiently small positive $\varepsilon_0(T)$ such that,
for $0<\varepsilon\leq\varepsilon_0(T)$, there is a solution $w$ of (16) such that
$w \in C^{\infty}([0,T]\times [1,R])$ and
$$\sup_{j+k\leq n}\Big\|
\Big(\frac{\partial}{\partial t}\Big)^j\Big(\frac{\partial}{\partial r}\Big)^k
w\Big\|_{L^{\infty}([0,T]\times[1,R])}\leq C_n\varepsilon, $$
or a solution $y \in C^{\infty}([0,T]\times[1,R])$ of (5) or (15) of the form
$$y(t,r)=\varepsilon y_1(t,r)+O(\varepsilon^2), $$
or a motion which can be expressed by the Lagrangian coordinate as
$$r(t,m)=\bar{r}(m)(1+\varepsilon y_1(t,\bar{r}(m))+O(\varepsilon^2))$$
for $0\leq t\leq T, 0\leq m\leq M$.
\end{Theorem}

Our task is to find the inverse image $\mathfrak{P}^{-1}(\varepsilon c)$
of the nonlinear mapping $\mathfrak{P}$ defined by
\begin{equation}
\mathfrak{P}(w)=\frac{\partial^2w}{\partial t^2}+
(1+\varepsilon a)\mathcal{L}w+\varepsilon b.
\end{equation}
Let us note that $\mathfrak{P}(0)=0$. This task, which will be done by applying the 
Nash-Moser theorem, will require a certain property of the derivative of $\mathfrak{P}$:
\begin{align}
D\mathfrak{P}(w)h&:=
\lim_{\tau\rightarrow 0}\frac{1}{\tau}(\mathfrak{P}(w+\tau h)-\mathfrak{P}(w)) \nonumber \\
&=\frac{\partial^2h}{\partial t^2}+(1+\varepsilon a_1)
\mathcal{L}h+
\varepsilon a_{21}r\frac{\partial h}{\partial r}+\varepsilon a_{20}h,
\end{align}
where
$$
a_1=a(t,r,w,r\frac{\partial w}{\partial r}, \varepsilon), \quad
a_{20}=\frac{\partial a}{\partial w}\mathcal{L}w
+\frac{\partial b}{\partial w}, \quad
a_{21}=\frac{\partial a}{\partial\Omega}\mathcal{L}w+\frac{\partial b}{\partial \Omega}
$$
are smooth functions of $t,r,w,\displaystyle r\frac{\partial w}{\partial r}, \varepsilon$.
Here $\Omega$ is the dummy of $\displaystyle r\frac{\partial w}{\partial r}$,
that is, $a$ and $b$ are functions of $t,r,w,\displaystyle\Omega=r\frac{\partial w}{\partial r},
 \varepsilon$ and
$\partial a/\partial\Omega \quad [(\partial b/\partial\Omega)]$
 denotes the partial derivative of $a \quad [(b)]$
with respect to $\Omega=r\partial w/\partial r$, respectively.
We consider $D\mathfrak{P}(w)$ as a second order linear
partial differential operator for each fixed $w$.

The following observation will play a crucial role in energy estimates later.\\

\textbullet {\it
We have
$$a_{21}=\frac{\gamma P}{\rho}
(1+y)^{-2\gamma+2}(1+y+v)^{-\gamma-2}\Big(
(\gamma+1)\frac{\partial^2Y}{\partial r^2}+
\frac{4\gamma}{r}\frac{\partial Y}{\partial r}+
\frac{2\varepsilon (\gamma-1)}{1+y}\Big(\frac{\partial Y}{\partial r}\Big)^2\Big),
$$
where
$$y=\varepsilon(y_1+w),\qquad v=r\frac{\partial y}{\partial r},
\qquad Y=y_1+w.
$$}

{\bf Proof}
It is easy to see
\begin{align*}
a_{21}&=(\partial_vG_I)\mathcal{L}Y+\varepsilon^{-1}\partial_vG_{II} \\
&=(\partial_vG_I)\Big(
-\frac{\gamma P}{\rho r}(3Y+V)'\Big)+\varepsilon^{-1}\frac{P}{\rho r^2}\partial_vG_{II0}+
\frac{1}{\gamma-1}\frac{1}{r^3}[U],
\end{align*}
where
$$[U]=\gamma(\partial_vG_I)(3Y+V)+
\partial_vG_I(-4Y)+
\varepsilon^{-1}\partial_vG_{II1},
$$
while $V$ stands for $\displaystyle r\frac{\partial Y}{\partial r}$.
Using
$$\partial_vG_I=(1+y)^2\frac{1}{\gamma}\partial_v^2G_2, \quad
\partial_vG_{II1}=(1+y)^2\frac{1}{\gamma}\partial_v^2G_2
((4-3\gamma)Y-\gamma V),$$
we see that $[U]=0$. Then a direct calculation leads us to the conclusion.
\hfill$\blacksquare$\\

 Hereafter we use the variable $x$ defined by (12) instead of $r$. We note
that
$$
x=R^2(R-r)+[R-r]_2, \quad
\frac{\partial}{\partial r}=-R^2(1+[x]_1)\frac{\partial}{\partial x}.
$$
Therefore a function of $1\leq r\leq R$ which is infinitely many times
continuously differentiable is also so as a function of $0\leq x\leq x_R$.

The consequence of the above observation is as follows.

We note that
$$\frac{\gamma P}{\rho }=\frac{1}{r}-\frac{1}{R}=\frac{x}{R^3}(1+[x]_1).$$
Therefore it follows from the above 
observation that there exists a smooth function $\hat{a}$ such that
$$\varepsilon a_{21}r\frac{\partial}{\partial r}=\varepsilon\hat{a}\cdot x\frac{\partial}{\partial x}.$$
Let us put
$$
b_1:=(1+\varepsilon a_1)L_1(x)+\varepsilon\hat{a}, \quad
b_0:=(1+\varepsilon a_1)L_0(x)+\varepsilon a_{20},
$$
taking into account the observation in Section 2, (13). Then we have

\begin{Lemma}
There are smooth functions $b_1,b_0$ of 
$t,x,w, \partial w/\partial x, \partial^2w/\partial x^2$ such that
$$D\mathfrak{P}(w)h=
\frac{\partial^2h}{\partial t^2}-(1+\varepsilon a_1)\triangle h
+b_1x\frac{\partial h}{\partial x}+b_0h.
$$
\end{Lemma}

{\bf Remark}\  The factor $x$ in the term $\displaystyle b_1x\frac{\partial h}{\partial x}$ is
important. In fact $\displaystyle B\frac{\partial h}{\partial x}$, $B$ being a non-zero constant,
without the factor
$x$ cannot be considered as a perturbation term, since
it has the same order with the principal part $\displaystyle \triangle h=x\frac{\partial^2h}{\partial x^2}+
\frac{N}{2}\frac{\partial h}{\partial x}$. See the proof of the following Lemma 3.\hfill$\square$

Using this representation of $D\mathfrak{P}(w)$, we can prove the following
energy estimate:

\begin{Lemma} 
If
a solution of
$D\mathfrak{P}(w)h=g$ satisfies 
$$h|_{x=x_R}=0, \quad h|_{t=0}=\frac{\partial h}{\partial t}\Big|_{t=0}=0,$$
then $h$ enjoys the energy inequality
$$\|\partial_th\|_{\mathfrak{X}}+\|\dot{D}h\|_{\mathfrak{X}}+\|h\|_{\mathfrak{X}}\leq
C \int_0^T\|g(t')\|_{\mathfrak{X}}dt',
$$
where $\partial_t=\partial/\partial t, \dot{D}=\sqrt{x}\partial/\partial x$ 
and $C$ depends only on $N, R, T$,
$A:=\|\varepsilon \partial_ta_1\|_{L^{\infty}}+\sqrt{2}\|\varepsilon \dot{D}a_1+b_1\|_{L^{\infty}}$
and
$B:=\|b_0\|_{L^{\infty}}$, provided that $|\varepsilon a_1|\leq 1/2$.
Here we have used the notation
$$\|y\|_{\mathfrak{X}}:=\Big(\int_0^{x_R}|y|^2x^{\frac{N}{2}-1}dx\Big)^{1/2}.$$
\end{Lemma}

{\bf Proof}
We consider the energy
$$E(t):=\int_0^{x_R}
((\partial_th)^2+(1+\varepsilon a_1)(\dot{D}h)^2)x^{\frac{N}{2}-1}dx.$$
Mutiplying the equation by $\partial_th$, and integrating by parts under the
boundary condition, we
get
\begin{align*}
\frac{1}{2}\frac{dE}{dt}&=
\int_0^{x_R}\Big(\frac{1}{2}\partial_t(\varepsilon a_1)(\dot{D}h)^2-
\dot{D}(\varepsilon a_1)(\dot{D}h)(\partial_th) + \\
&-\sqrt{x}b_1(\dot{D}h)(\partial_th)-
b_0h(\partial_th)+
g(\partial_th)\Big)x^{\frac{N}{2}-1}dx, 
\end{align*}
which implies
$$\frac{1}{2}\frac{dE}{dt}\leq
AE+B\Big|\int_0^{x_R}h(\partial_th)x^{\frac{N}{2}-1}dx\Big|+
E^{1/2}\|g(t)\|_{\mathfrak{X}}.$$
On the other hand, using the initial condition, we see that
$U(t):=\|h\|_{\mathfrak{X}}^2$ enjoys
$$
\frac{1}{2}\frac{dU}{dt}=\int_0^{x_R} h(\partial_th)x^{\frac{N}{2}-1}dx \leq U^{1/2}E^{1/2}, \quad
U(0)=0.
$$
Hence $U(t)\leq \int_0^tE^{1/2}$ and
$$\Big|\int_0^{x_R} h(\partial_th)x^{\frac{N}{2}-1}dx
\Big|\leq E^{1/2}(t)\int_0^tE^{1/2}.$$
Summing up, we have
$$\frac{1}{2}\frac{dE}{dt}\leq
AE+BE(t)^{1/2}\int_0^tE^{1/2}+E^{1/2}\|g(t)\|_{\mathfrak{X}},\quad
E(0)=0.
$$
By the Gronwall's lemma, we can derive the inequality
$$E^{1/2}(t)\leq C\int_0^t\|g(t')\|_{\mathfrak{X}}dt'.$$
\hfill$\blacksquare$\\

\section{Proof of the main result}

We are going to apply the Nash-Moser theorem formulated by R. Hamilton 
(\cite[p.171, III.1.1.1]{Hamilton}):\\

\noindent{\bf Nash-Moser Theorem}\  {\it Let $\mathfrak{E}_0$ and $\mathfrak{E}$ be tame spaces, $U$ an open subset of 
$\mathfrak{E}_0$ and $\mathfrak{P}: U\rightarrow \mathfrak{E}$ a smooth tame map. 
Suppose that the equation for the derivative
$D\mathfrak{P}(w)h=g$ has a unique solution 
$h=V\mathfrak{P}(w)g$ for all $w$ in $U$ and all $g$, and
that the family of inverse $V\mathfrak{P}: U\times \mathfrak{E}
\rightarrow \mathfrak{E}_0$ is a smooth tame map. Then $\mathfrak{P}$ is
locally invertible.}\\

Let us recall the definitions of `tame spaces' and `tame maps' in the sense of \cite{Hamilton}:

\begin{df}

1) A {\bf graded space} $\mathfrak{E}$ is a 
Fr\'{e}chet space whose topology is given
by a grading of
seminorms $(\|\cdot\|_{\mathfrak{E},n})_{n\in\mathbb{N}}$ such that 
$\|y\|_{\mathfrak{E},n}\leq \|y\|_{\mathfrak{E},n+1}$;

2) A linear map $L$ from a graded space $\mathfrak{E}$ into a graded space $\mathfrak{F}$
is said to be {\bf tame} if $\|Ly\|_{\mathfrak{F},n}
\leq C_n\|y\|_{\mathfrak{E}, n+r}$ for all $n$ with some $r$;

3) A graded space $\mathfrak{E}$ is said to be a {\bf tame direct summand} of a graded
space $\mathfrak{F}$ if we can find tame linear maps $L:\mathfrak{E}\rightarrow
\mathfrak{F}$ and $M:\mathfrak{F}\rightarrow \mathfrak{E}$ such that the composition
$M\circ L$ is the identity of $\mathfrak{E}$;

4) A graded space $\mathfrak{E}$ is said to be {\bf tame} if it is a tame direct summand 
of $\Sigma(B)$, where $B$ is a Banach space and
$$
\Sigma(B)=\{ f=(f_k)_{k\in\mathbb{N}}\in B^{\mathbb{N}}\quad | \quad
\|f\|_{\Sigma(B),n}:=\sum_ke^{nk}\|f_k\|_B <\infty \quad \forall n\};
$$

5) A continuous mapping $\mathfrak{P}: U\rightarrow \mathfrak{F}$, where
$U$ is an open subset of $\mathfrak{E}$, $\mathfrak{E},\mathfrak{F}$ being graded spaces,
is said to be {\bf tame} if $\mathfrak{P}$ satisfies a tame estimate, that is,
$$\|\mathfrak{P}(w)\|_{\mathfrak{F},n}\leq C_n(1+\|w\|_{\mathfrak{E},n+r})$$
for all $n$ with some $r$.
\end{df}

This section is devoted to set a framework to apply the above Nash-Moser theorem,
say,
a tame map $\mathfrak{P}$ from an open set of a tame
space $\mathfrak{E}_0$ into a tame space $\mathfrak{E}$. A tame
estimate
of $V\mathfrak{P}$ will be verified in the next section.

In order to apply the Nash-Moser theorem, we consider the spaces of
functions of $t$ and $x$:
\begin{align*}
\mathfrak{E}&:= \{ y\in C^{\infty}([-2\tau_1,T]\times[0,x_R])\quad | \quad y(t,x)=0 \quad\mbox{for}\quad
-2\tau_1\leq t\leq -\tau_1 \}, \\
\mathfrak{E}_0&:=\{ w \in \mathfrak{E}\ |\quad w|_{x=x_R}=0\}.
\end{align*}
Here $\tau_1$ is a positive number.
Let $U$ be the set of all functions $w$ in $\mathfrak{E}_0$ such that
$|w|+|\partial w/\partial x| <1$ and suppose that $|\varepsilon|\leq
\varepsilon_1$, $\varepsilon_1$ being a small positive number. Then
we can consider that the nonlinear mapping $\mathfrak{P}$ maps $U(\subset\mathfrak{E}_0)$
into $\mathfrak{E}$, since the coefficients $a,b$ are smooth functions of $t,
x, \varepsilon w$ and $\varepsilon\partial w/\partial x$. 
Let us assume that $\varepsilon c(t,x)=0$ for 
$-2\tau_1\leq t\leq -\tau_1$ after changing the value of $c$ 
for $-2\tau_1\leq t <0$. To fix the idea, we replace $c(t,x)$
by $\alpha(t)c(t,x)$ with a cut off function
$\alpha\in C^{\infty}(\mathbb{R})$ such that
$\alpha(t)=1$ for $t\geq 0$ and $\alpha(t)=0$ for $t \leq -\tau_1$.
Then $\mathfrak{P}^{-1}(\varepsilon c)$ is
a solution of (16) on $t\geq 0$.\\

We should show that the Fr\'{e}chet space $\mathfrak{E}$ is tame for some gradings of norms.
For $y\in \mathfrak{E}$, $n\in\mathbb{N}$, let us define
$$
\|y\|_n^{(\infty)}:=
\sup_{0\le j+k\leq n}\Big\|\Big(-\frac{\partial^2}{\partial t^2}
\Big)^{j}(-\triangle)^{k}y\Big\|_{L^{\infty}([-2\tau_1,T]\times[0,x_R])}.
$$
Then we can claim that $\mathfrak{E}$ turns out to be tame by this
grading $(\|\cdot\|_{n}^{(\infty)})_n$ (see \cite[p.136, II.1.3.6 and p.137, II 1.3.7]{Hamilton}).
In fact, even if $N$ is not an integer, we can define the Fourier transformation $Fy(\xi)$ of a function $y(x)$
 for $0\leq x<\infty$ by
$$Fy(\xi):=\int_0^{\infty}K(\xi x)y(x)x^{\frac{N}{2}-1}dx. $$
Here $K(X)$ is an entire function of $X \in \mathbb{C}$ given by
$$K(X)=2(\sqrt{X})^{-\frac{N}{2}+1}J_{\frac{N}{2}-1}(4\sqrt{X})
=2^{\frac{N}{2}}
\sum_{k=0}^{\infty}\frac{(-4)^kX^k}{k!\Gamma(\frac{N}{2}+k)}, $$
$J_{\frac{N}{2}-1}$ being the Bessel function of order $\frac{N}{2}-1$. Then we have
$$F({-\triangle y})(\xi)=4\xi\cdot F{y}(\xi) $$
and the inverse of the transformation $F$ is $F$ itself. See \cite[p.65]{Sneddon}. 
Then it is easy to see $\mathfrak{E}$ endowed with
the grading $(\|y\|_n^{(\infty)})_n$ 
is a tame direct summand of the tame space
$$\mathfrak{F}:=L_1^{\infty}(\mathbb{R}\times [0,\infty), d\tau\otimes 
\xi^{\frac{N}{2}-1}d\xi, \log (1+\tau^2+4\xi)) $$
through the Fourier transformation
$$\mathcal{F}y(\tau, \xi)=
\frac{1}{\sqrt{2\pi}}\int e^{-\sqrt{-1}\tau t}Fy(t,\cdot)(\xi)dt$$
and its inverse
applied to the space $\tilde{\mathfrak{E}}_0:=C_0^{\infty}((-2T-2\tau_1,2T)\times [0, x_R+1))$, 
into which functions of
$\mathfrak{E}$ can be extended (see, e.g. 
\cite[p.88,  4.28 Theorem]{Adams}, -the existence of `total extension operator') and the space
$$\dot{\mathfrak{E}}:=\dot{C}^{\infty}(\mathbb{R}\times [0,\infty)) :=
\{y | \forall j\quad\forall k \quad\lim_{L\rightarrow\infty}\sup_{|t|\geq L,x\geq L}|(-\partial_t^2)^j(-\triangle)^ky|=0\},$$
for which functions of $\mathfrak{E}$ are restrictions.
Actually, if we denote by $\mathfrak{e}:\mathfrak{E}\rightarrow
\tilde{\mathfrak{E}}_0$ the extension operator, and by
$\mathfrak{r}:\dot{\mathfrak{E}}\rightarrow\mathfrak{E}$ the restriction operator,
then
the operators $\mathcal{F}\circ\mathfrak{e}:\mathfrak{E}\rightarrow \mathfrak{F}$ and
$\mathfrak{r}\circ\mathcal{F}:\mathfrak{F}\rightarrow \mathfrak{E}$ are tame and
the composition $(\mathfrak{r}\circ\mathcal{F})\circ(\mathcal{F}\circ\mathfrak{e})$ is
the identity of $\mathfrak{E}$. 
For the details, see the proof of \cite[p.137, II.1.3.6.Theorem]{Hamilton}.
This shows that $\mathfrak{E}$ is tame with respect to
the grading $(\|\cdot\|_n^{(\infty)})_n$.\\

On the other hand, let us define
$$\|y\|_n^{(2)}:=\Big(
\sum_{0\le j+k\leq n}
\int_{-\tau_1}^T\|\Big(-\frac{\partial^2}{\partial t^2}\Big)^{j}(-\triangle)^{k}y\|_{\mathfrak{X}}^2dt
\Big)^{1/2}. $$
Here $\mathfrak{X}=L^2((0,x_R); x^{\frac{N}{2}-1}dx)$ and
$$\|y\|_{\mathfrak{X}}:=\Big(\int_0^{x_R}|y(x)|^2x^{\frac{N}{2}-1}dx\Big)^{1/2}.$$
We have
$$\sqrt{\frac{N}{2}}\|y\|_{\mathfrak{X}}\leq \|y\|_{L^{\infty}}\leq C
\sup_{j\leq\sigma}\|(-\triangle)^jy\|_{\mathfrak{X}},$$
by the Sobolev imbedding theorem (see Appendix A), provided that $2\sigma > N/2$. 
The derivatives with respect to $t$
can be treated more simply. Then
we see that
the grading $(\|\cdot\|_n^{(2)})_n$ is tamely equivalent to
the grading $(\|\cdot\|_n^{(\infty)})_n$,
that is, we have
$$\frac{1}{C}\|y\|_n^{(2)}\leq\|y\|_n^{(\infty)}\leq
C\|y\|_{n+s}^{(2)}$$
with $2s >1+N/2$. Hence
$\mathfrak{E}$ is tame with respect to $(\|\cdot\|_n^{(2)})_n$, too.
The grading $(\|\cdot\|_n^{(2)})_n$ will be suitable for estimates
of solutions of the associated linear wave equations.

Note that $\mathfrak{E}_0$ is a closed subspace of $\mathfrak{E}$ endowed with these
gradings.
\medskip

 Now we verify the nonlinear mapping $\mathfrak{P}$
is tame for the grading
$(\|\cdot\|_n^{(\infty)})_n$. To do so,
we write
$$\mathfrak{P}(w)=F(t, x, Dw, \triangle w, w_{tt}), $$
where $D=\partial/\partial x$, $F$ is a smooth function of
$t,x, Dw, \triangle w, w_{tt}$ and linear in $\triangle w, w_{tt}$.
According to \cite{Hamilton} (see p.142, II.2.1.6 and p.145, II.2.2.6), 
it is sufficient to prove the
linear differential operator $w \mapsto Dw=\partial w/\partial x$ is tame.
But it is clear because of the following result.
\begin{Proposition}
For any $m\in\mathbb{N}$ and for any $y\in C^{\infty}[0,1]$ we have the formula
$$\triangle^mDy(x)=
x^{-\frac{N}{2}-m-1}\int_0^x
\triangle^{m+1}y(x')(x')^{\frac{N}{2}+m}dx'. $$
As a corollary it holds that, for any $m,k \in \mathbb{N}$,
$$\|(-\triangle)^mD^ky\|_{L^{\infty}}\leq
\frac{1}{\prod_{j=0}^{k-1}(\frac{N}{2}+m+j)}\|(-\triangle)^{m+k}y\|_{L^{\infty}}.$$
\end{Proposition}
{\bf Proof} It is easy by integration by parts in induction on $m$ starting from the formula
$$Dy(x)=
x^{-\frac{N}{2}}\int_0^x
\triangle y(x')(x')^{\frac{N}{2}-1}dx'. $$
\hfill$\blacksquare$
\\

In parallel with the results of \cite{Hamilton} (see p.144, II.2.2.3.Corollary 
and p.145, II.2.2.5.Theorem), we
should use the following two propositions. Proofs for these propositions are given in Appendix B.

\begin{Proposition}
For any positive integer $m$, there is a constant $C$ such that
$$|\triangle^m(f\cdot g)|_0\leq
C(|\triangle^mf|_0|g|_0+
|f|_0|\triangle^mg|_0), $$
where $|\cdot |_0$ stands for $\|\cdot \|_{L^{\infty}}$.
\end{Proposition}

\begin{Proposition}
Let $F(x,y)$ be a smooth function of $x$ and $y$ and $M$ be a positive number. Then
for any positive integer $m$, there is a constant $C>0$ such that
$$|\triangle^mF(x,y(x))|_0\leq C
(1+|y|_m^*) $$
provided that $|y|_0\leq M$, where we denote
$$|y|_m^*=\sup_{0\leq j\leq m}\|(-\triangle)^jy\|_{L^{\infty}}.$$
\end{Proposition}

Summing up, we can claim that
$$\|\mathfrak{P}(w)\|_n^{(\infty)}\leq C(1+\|w\|_{n+1}^{(\infty)}),$$
provided that $\|w\|_1^{(\infty)}\leq M$.
This says that the mapping $\mathfrak{P}$ is tame with respect to the grading
$(\|\cdot\|_n^{(\infty)})_n$.\\

Therefore the problem is concentrated to estimates of 
the solution and its higher derivatives of the linear equation
$$D\mathfrak{P}(w)h=g,$$
when $w$ is fixed in $\mathfrak{E}_0$ and $g$ is given in $\mathfrak{E}$.
A tame estimate of the mapping $(w,g)\mapsto h$ will be discussed in the next section. 
This will completes the proof of the main result.

\section{Tame estimate of solutions of linear wave equations}

This section is devoted to verify a tame estimate
$$\|h\|_n^{(2)}\leq C(1+\|g\|_n^{(2)}+\|w\|_{n+3+s}^{(2)})$$
with $2s>1+N/2$, provided that
$\|g\|_1^{(2)}\leq M$ and
$\|w\|_{3+s}^{(2)}\leq M$. Here $h$ is the solution of the
equation
$$D\mathfrak{P}(w)h\equiv
\frac{\partial^2h}{\partial t^2}-(1+\varepsilon a_1)\triangle h+
b_1x\frac{\partial h}{\partial x}+b_0h=g$$
for given $g\in\mathfrak{E}$, provided that $|\varepsilon a_1|\leq 1/2$.
This estimate says that the mapping
$(w,g)\mapsto h$ is tame with respect to the grading
$(\|\cdot\|_n^{(2)})_n$.

Therefore we are considering the wave equation
\begin{equation}
\frac{\partial^2h}{\partial t^2}+\mathcal{A}h=g(t,x), \qquad (0\leq t\leq T, 0\leq x\leq 1),
\end{equation}
where
$$
\mathcal{A}h=-b_2\triangle h +b_1\check{D}h+b_0h, \quad
\triangle =x\frac{d^2}{dx^2}+\frac{N}{2}\frac{d}{dx} \quad \check{D}=x\frac{d}{dx}.
$$
We denote $\vec{b}=(b_2,b_1,b_0)$. The given function $\vec{b}(t,x)$ is supposed to be in
$C^{\infty}([0,T]\times[0,1])$ and we assume that
$|b_2(t,x)-1|\leq 1/2$. The function $g(t,x)$
belongs to $C^{\infty}([0,T]\times[0,1])$ and
we suppose that
\begin{equation}
g(t,x)=0 \qquad \mbox{for}\qquad 0\leq t\leq \tau_1,
\end{equation}
where $\tau_1$ is a positive number. 

In this section the $x$-interval $[0,x_R]$ has been 
normalized as
$[0,1]$ without loss of generality, and the parallel
translation $t \rightarrow t+2\tau_1$ has been done.

Let us consider the
initial boundary value problem (IBP):
$$
\frac{\partial^2h}{\partial t^2}+\mathcal{A}h=g(t,x), \quad
h|_{x=1}=0, \quad
h|_{t=0}=\frac{\partial h}{\partial t}\Big|_{t=0}=0 .
$$

Then (IBP) admits a unique solution $h(t,x)$ thanks to the energy estimate, and
$h(t,x)=0$ for $0\leq t\leq \tau_1$ because of the uniqueness.
Moreover, since the compatibility conditions are satisfied, the unique
solution turns out to be smooth. A proof can be found e.g. in
\cite[Chapter 2]{Ikawa}, at least for the case in which the coefficients of
$\mathcal{A}$ do not depend on $t$.
To satisfy ourselves, we shall give a brief sketch 
of a proof of the existence of smooth solutions in Appendix C.
We are going to get estimates of the higher
derivatives of $h$ by them of $g$ and the coefficients $b_2,b_1,b_0$.\\

\subsection{Notations}

It may be difficult to deduce the required tame estimate by using the norms
$\|\cdot\|_n^{(2)}$ directly. Hence let us introduce
auxiliary other norms
$\|\cdot\|_n, |\cdot|_n$ and $\|\cdot\|_n^T, |\cdot|_n^T$ defined as follows.

For $m,n\in\mathbb{N}$ and for functions $y=y(x)$ of
$x \in [0,1]$, we put
\begin{align*}
\langle y\rangle_{2m} &:=\|\triangle^my\|, 
\qquad \|y\|:=\|y\|_{\mathfrak{X}}:=\Big(\int_0^1|y(x)|^2x^{\frac{N}{2}-1}dx\Big)^{1/2}, \\
\langle y\rangle_{2m+1} &:=\|\dot{D}\triangle^my\|, \quad \dot{D}=\sqrt{x}\frac{d}{dx}, \quad
\|y\|_n:=\Big(\sum_{0\leq\ell\leq n}\langle y\rangle_{\ell}^2\Big)^{1/2}, \\
|y|_n&:=\max_{0\leq\ell\leq n}
\|\dot{D}^{\ell}y\|_{L^{\infty}(0,1)}.
\end{align*}

For $n\in\mathbb{N}$, a fixed $T>0$, and for functions
$y=y(t,x)$ of $(t,x)\in[0,T]\times[0,1]$, we put
$$
\|y\|_n^T:=\Big(\sum_{j+k\leq n}
\int_0^T\|\partial_t^jy\|_k^2dt\Big)^{1/2}, \quad
|y|_n^T:=\max_{j+k\leq n}\|\partial_t^j\dot{D}^ky\|_{L^{\infty}([0,T]\times[0,1])}.
$$
Here $\partial_t=\partial/\partial t$. 

\begin{df}Let us say that a grading of norms $(p_n)_{n\in\mathbb{N}}$ is {\bf interpolation admissible} if
for $\ell\leq m\leq n$ it holds that
$$p_m(f)\leq Cp_n(f)^{\frac{m-\ell}{n-\ell}}p_{\ell}(f)^{\frac{n-m}{n-\ell}}.
$$
\end{df}

It is well known that, if and only if
$$p_n(f)^2\leq C p_{n+1}(f)p_{n-1}(f) $$
for any $n\geq 1$, $(p_n)_n$ is interpolation admissible.
If $(p_n)_n$ and $(q_n)_n$ are interpolation admissible, and if
$(i,j)$ lies on the line segment joining
$(k,\ell)$ and $(m,n)$, then
$$p_i(f)q_j(g)\leq C
(p_k(f)q_{\ell}(g)+p_m(f)q_n(g)).
$$ (For a proof, see \cite[ p.144, 2.2.2. Corollary]{Hamilton}.)

It is well-known that $(|\cdot|_n)_n$ and
$(|\cdot|_n^T)_n$ are
interpolation admissible, since $\dot{D}=\partial/\partial\xi$, where
$x=\xi^2/4$.

Moreover $(\|\cdot\|_n)_n$ and $(\|\cdot\|_n^T)_n$ are interpolation admissible. To verify it, it is
sufficient to note that
$y=\sum_{k=1}^{\infty}c_k\phi_k \in C_0^{\infty}[0,1)$
enjoys
$\langle y\rangle_{\ell}=\Big(\sum_k\lambda_k^{\ell}|c_k|^2\Big)^{1/2}. $
Here $(\lambda_k)_k$ are eigenvalues of $-\triangle$ with the Dirichlet
boundary condition at $x=1$ and $(\phi_k)_k$ are associated eigenfunctions.
We note that $(\dot{D}\phi_n/\sqrt{\lambda_n})_{n=1,2,...}$ is a complete
orthonormal system of $\mathfrak{X}$
and $(\dot{D}y|\dot{D}\phi)_{\mathfrak{X}}=
(-\triangle y|\phi)_{\mathfrak{X}}$ if $y \in C^{\infty}[0,1)$.

 Then it is clear by the Schwartz inequality that
$$\langle y\rangle_n^2\leq \langle y\rangle_{n+1}\langle y\rangle_{n-1}$$
for $y \in C_0^{\infty}[0,1)$. Since
$\langle y\rangle_j\leq \langle y\rangle_{j'}$ for $j\leq j', y\in C_0^{\infty}[0,1)$, we have
$$\langle y\rangle_{\ell}\leq \|y\|_{\ell}\leq C\cdot \langle y\rangle_{\ell}$$
and
$$\|y\|_n^2\leq C\|y\|_{n-1}\|y\|_{n+1}$$
at least for $y \in C_0^{\infty}[0,1)$.
By using a continuous linear extension of functions 
on $[0,1]$ to functions on $[0,2]$ with supports in $[0,3/2)$, we can
claim that this inequality holds for any $y \in C^{\infty}[0,1]$
with a suitable change of the constant $C$. We refer to
\cite[Chapter 3, Section 4, Theorem 3.11]{Mizohata}. It is sufficient to note
the following
\begin{Proposition}
If $\alpha(x) \in C^{\infty}(\mathbb{R})$ is fixed, then
there is a constant $C$ depending on $\alpha$ such that 
$$\|\alpha y\|_n\leq C \|y\|_n.$$
\end{Proposition}

A proof can be found in Appendix B. Hence $(\|\cdot\|_n)_n$
and $(\|\cdot\|_n^T)_n$ are interpolation admissible.\\

\subsection{Goal of this Section}

Using the norms $\|\cdot\|_n^T$ and $|\cdot|_n^T$, we state our goal of
this section as the following

\begin{Lemma}
 Assume that $|b_2-1|\leq 1/2$, 
$|\vec{b}|_2^T\leq M$
and $\|g\|_1^T \leq M$. Then there is a constant $C_n=C_n(T,M,N)$ such that
if $h$ is the solution of (IBP) then 
$$\|h\|_{n+2}^T\leq C_n
(1+\|g\|_{n+1}^T+|\vec{b}|_{n+3}^T).$$
\end{Lemma}

Let us note that this lemma implies the required tame estimate.
We see that $\|y\|_{2m}^T$ is equivalent to
$$\|y\|_m^{(2)}=\Big(\sum_{j+k\leq m}
\int_0^T\langle\partial_t^{2j}y\rangle_{2k}^2dt\Big)^{1/2}$$
for $y \in C^{\infty}([0,T]\times[0,1])$. In fact 
it is sufficient to note the following

\begin{Proposition}
For any $y \in C^{\infty}([0,1])$ we have
$$\|\dot{D}\triangle^my\|_{\mathfrak{X}}\leq C(\|\triangle^my\|_{\mathfrak{X}}+
\|\triangle^{m+1}y\|_{\mathfrak{X}}).$$
\end{Proposition}

A proof can be found in Appendix B. 
Therefore the conclusion of the Lemma reads:
$$\|h\|_{m}^{(2)}\leq C(1+\|g\|_{m}^{(2)}+\|w\|_{m+3+s}^{(2)})$$
with $2s>1+N/2$,
provided that $\|g\|_1^{(2)}\leq M$ and $\|w\|_{3+s}^{(2)}\leq M$,
since $\vec{b}$ is a smooth function of
$w, Dw, D^2w, \partial_t^2w$ in our context
so that $|\vec{b}|_{n+3}^T\leq C(1+|w|_{n+7}^T)$,
provided that $|w|_4^T\leq M'$. 
Note that $|w|_{n+7}^T\leq C\|w\|_{m+3}^{(\infty)}$ if $2m=n+2$. This is
the required tame estimate.\\

\subsection{Elliptic a priori estimates}

In order to prove Lemma 4, we shall use the so called elliptic
a priori estimate of the operator $\mathcal{A}$ of G\aa rding's type:

\begin{Proposition}
Suppose $|b_2-1|\leq 1/2$ and $|\vec{b}|_2\leq M$. Then
$$\|y\|_{n+2}
\leq C(\|\mathcal{A}y\|_n+\|y\|_1+|\vec{b}|_{n+3}\|y\|).$$
\end{Proposition}

A proof of this proposition can 
be found in Appendix D.

\subsection{Proof of Lemma 4}

Let us go to prove Lemma 4, using Proposition 8 and the energy estimate (Lemma 3).

The essence of the proof lies on the fact that
all higher derivatives $\partial_t^jh$ of the solution $h$ with respect to $t$
satisfy the Dirichlet boundary condition, and therefore enjoy the energy estimate. The apparent 
complicatedness of the discussion comes from
the situation that the coefficients of $\mathcal{A}$ depend on $t$. The logical structure
of the discussion would be quite simple and clear if the coefficients 
were constants with respect to $t$. See, e.g., the discussion in \cite[\S 2.2 (c)]{Ikawa}.

Hereafter we generally denote by
$H$ a solution of the boundary
value problem
$$\frac{\partial^2H}{\partial t^2}+\mathcal{A}H=G(t,x),\qquad H|_{x=1}=0 $$
such that $H(t,x)=0$ for $0\leq t\leq \tau_1$. Thus $\partial_t^jH|_{t=0}=0$ for
any $j\in\mathbb{N}$. The time derivative $H_j=\partial_t^jH$ satisfies
$$\frac{\partial^2H_j}{\partial t^2} +\mathcal{A}H_j=G_j,\qquad H_j|_{x=1}=0,$$
where
$$G_j:=\partial_t^jG-[\partial_t^j,\mathcal{A}]H.$$
We put $G_0=G$. Hereafter we always assume that $|b_2-1|\leq 1/2$ and $|\vec{b}|_2^T\leq M$.

Note that we have the energy
estimate
$$\|\partial_tH\|+\|H\|_1\leq C\int_0^t\|G(t')\|dt'$$
for $0\leq t\leq T$. (See Lemma 3.)

{\bf Remark} In this subsection $H$ and $G$ do not mean the
particular functions defined in Section 1.\hfill$\square$

First let us reduce the estimates of the mixed derivatives $\|\partial_t^jH\|_k$
to those of purely time derivatives $\|\partial_t^jH\|_1$ by
using the elliptic a priori estimate and the equation. In fact,
putting
$$Z_n(H):=\sum_{j+k=n}\|\partial_t^jH\|_k,$$
we can claim the following

\begin{Proposition}
For $n\in\mathbb{N}$ we have
$$
Z_{n+2}(H)\leq C
(\|\partial_t^{n+1}H\|_1+
\sum_{j+k=n}\|G_j\|_k+\sum_{j+k=n}(\|\partial_t^jH\|_1+|\vec{b}|_{k+3}\|\partial_t^jH\|).
$$
\end{Proposition}

{\bf Proof}\ First we show that
for $n\in\mathbb{N}$ we have
$$Z_{n+2}(H)
\leq C(Z_{n+1}(\partial_tH)+\|G\|_n+
\|H\|_1+|\vec{b}|_{n+3}\|H\|).$$

In fact by definition we have
$$Z_{n+2}(H)=Z_{n+1}(\partial_tH)+\|H\|_{n+2}.$$
By Proposition 8 we have
\begin{align*}
\|H\|_{n+2}&\leq C(\|\mathcal{A}H\|_n+\|H\|_1+|\vec{b}|_{n+3}\|H\|) \\
&\leq C(\|\partial_t^2H-G\|_n+\|H\|_1+|\vec{b}|_{n+3}\|H\|) \\
&\leq C(\|\partial_t^2H\|_n+\|G\|_n+\|H\|_1+|\vec{b}|_{n+3}\|H\|).
\end{align*}
Note that $\|\partial_t^2H\|_{n}\leq Z_{n+1}(\partial_tH)$. 

This implies by induction the desired estimates. \hfill$\blacksquare$\\

Let us apply the energy estimate to $\|\partial_t^{n+1}H\|_1$. Then we get

\begin{Proposition} We have
\begin{align}
\|H\|_{n+2}^T&\leq C\Big(\|G\|_{n+1}^T+
\sum_{j+k\leq n}\Big(\int_0^T\|G_j\|_k^2\Big)^{1/2}+ \nonumber \\
&+\sum_{0\leq j\leq n}\sup_{0\leq t\leq T}\|\partial_t^jH\|_1+\|H\|_{n+1}^T+
|\vec{b}|_{n+3}^T\|H\|^T\Big).
\end{align}
\end{Proposition}

{\bf Proof}\ The energy estimate of $\partial_t^{n+1}H$ reads
$$\|\partial_t^{n+1}H\|_1\leq C
\Big(\int_0^t\|\partial_t^{n+1}G\|+\int_0^t
\|[\partial_t^{n+1},\mathcal{A}]H\| \Big).
$$
But
\begin{align*}
\int_0^t\|[\partial_t^{n+1},\mathcal{A}]H\| &
\leq C\sum_{\alpha+\beta=n+1, \alpha\not=0,}|\partial_t^{\alpha}\vec{b}|_0^T
\Big(\int_0^t
\|\partial_t^{\beta}H\|_2^2\Big)^{1/2} \\
&\leq C'\Big(|\vec{b}|_1^T
\Big(\int_0^tZ_{n+2}(H)^2\Big)^{1/2}+|\vec{b}|_{n+1}^T\|H\|_2^T\Big)
\end{align*}
by interpolation. 
Then Proposition 9 implies
$$Z_{n+2}(H)(t)\leq
C\Big(\Big(\int_0^t
Z_{n+2}(H)^2\Big)^{1/2}
+F_n(t)\Big),$$
where
\begin{align*}
F_n(t)&=\int_0^t\|\partial_t^{n+1}G(t')\|dt'+
|\vec{b}|_{n+1}^T\|H\|_2^T+\sum_{j+k=n}\|G_j\|_k+\\
&+\sum_{j+k=n}(\|\partial_t^jH\|_1+|\vec{b}|_{k+3}\|\partial_t^jH\|).
\end{align*}
 
We can apply the Gronwall's lemma to this inequality. The result is
$$
Z_{n+2}(H)(t)\leq
C\Big(\Big(\int_0^tF_n(t')^2dt'\Big)^{1/2}+F_n(t)\Big).
$$

Integrating the above inequality, we see
\begin{align*}
\|H\|_{n+2}^T&=\Big(\sum_{j+k\leq n+2}\int_0^T\|\partial_t^jH\|_k^2dt\Big)^{1/2} \\
&= \Big(\int_0^T(\|H\|^2+\|\partial_tH\|^2+\|H\|_1^2+\sum_{0\leq\nu\leq n}Z_{\nu+2}(H)^2)dt\Big)^{1/2} \\
&\leq  \Big(\int_0^T(\|H\|^2+\|\partial_tH\|^2+\|H\|_1^2)dt+
C\sum_{0\leq\nu\leq n}\int_0^TF_{\nu}^2\Big)^{1/2} \\
&\leq C'\Big(\|G\|_{n+1}^T+|\vec{b}|_{n+1}^T\|H\|_2^T+
\sum_{j+k\leq n}\Big(\int_0^T\|G_j\|_k^2\Big)^{1/2}+  \\
&+\sum_{0\leq j\leq n}\sup_{0\leq t\leq T}\|\partial_t^jH\|_1+
|\vec{b}|_{2}^T\|H\|_{n+1}^T+
|\vec{b}|_{n+3}^T\|H\|^T\Big)
\end{align*}
by interpolation.
Hereafter we suppose that $n\geq 1$. Then by interpolation we have
$$|\vec{b}|_{n+1}^T\|H\|_2^T\leq C(|\vec{b}|_2^T\|H\|_{n+1}^T+
|\vec{b}|_{n+3}^T\|H\|^T)$$
and therefore this completes the proof of (21). \hfill $\blacksquare$\\

Now, let us estimate the second and third terms in the right-hand side of (21).
First we have
\begin{equation}
\sum_{j+k=\nu}\Big(\int_0^T\|G_j\|_k^2\Big)^{1/2}\leq C(\|G\|_{\nu}^T +\|H\|_{\nu+1}^T 
+|\vec{b}|_{\nu+3}^T\|H\|^T)
\end{equation}

{\bf Proof} It is sufficient to estimate $\displaystyle
\int_0^T\|[\partial_t^j, \mathcal{A}(\vec{b})]H\|_k^2dt$.
But
$$[\partial_t^j, \mathcal{A}(\vec{b})]H=
\sum_{\alpha+\beta=j,\alpha\not=0}
\binom{j}{\alpha}
\mathcal{A}(\partial_t^{\alpha}\vec{b})
\partial_t^{\beta}H, $$
and
$$\|\mathcal{A}(\partial_t^{\alpha}\vec{b})\partial_t^{\beta}H\|_k
\leq C(\|\partial_t^{\beta}H\|_{k+2}+
|\partial_t^{\alpha}\vec{b}|_{k+3}\|\partial_t^{\beta}H\|),$$
since
$$\|\mathcal{A}(\vec{b})y\|_k\leq C
(\|y\|_{k+2}+|\vec{b}|_{k+3}\|y\|).$$
(The estimate of
$\|\mathcal{A}y\|_n$ can be derived by the discussion of the preceding subsection,
keeping in mind that
$\triangle^m\mathcal{A}=\mathcal{A}\triangle^m+[\triangle^m,\mathcal{A}]$.)
By interpolation, we have, for $\alpha+\beta+k=\nu, \alpha\not=0$,
\begin{align*}
\Big(\int_0^T\|\mathcal{A}(\partial_t^{\alpha}\vec{b})
\partial_t^{\beta}H\|_k^2\Big)^{1/2}&\leq
C(\|H\|_{\beta+k+2}^T+|\vec{b}|_{\alpha+k+3}^T\|H\|_{\beta}^T) \\
&\leq C'(\|H\|_{\nu+1}^T+
|\vec{b}|_2^T\|H\|_{\nu+1}^T+
|\vec{b}|_{\nu+3}^T\|H\|^T).
\end{align*}
\hfill$\blacksquare$

Next we have
\begin{equation}
\sup_{0\leq t\leq T}\|\partial_t^jH\|_1 \leq
C(\|G\|_j^T + \|H\|_{j+1}^T +|\vec{b}|_{j+3}^T\|H\|^T).
\end{equation}

{\bf Proof} By the energy estimate, we have
$$\|\partial_t^jH\|_1\leq C\int_0^T\|G_j\|.$$
Here we can use the estimate of $\displaystyle
\int_0^T\|[\partial_t^j, \mathcal{A}]H\|^2$
given in the proof of the preceding proposition with
$k=0, n=j$. \hfill$\blacksquare$

Substituting (22),(23) to (21) for $h=H$, we have
$$
\|h\|_{n+2}^T\leq
C(\|g\|_{n+1}^T+\|h\|_{n+1}^T+|\vec{b}|_{n+3}^T\|h\|^T).
$$

Noting that 
$$\|h\|^T\leq C\|g\|^T\leq CM,$$
we have
$$
\|h\|_{n+2}^T\leq C
(\|h\|_{n+1}^T+\|g\|_{n+1}^T+|\vec{b}|_{n+3}^T),
$$ 
which implies inductively that
$$
\|h\|_{n+2}^T\leq C
(1+\|g\|_{n+1}^T+|\vec{b}|_{n+3}^T),
$$
provided that $\|g\|_1^T\leq M$ and $|\vec{b}|_2^T\leq M$. 

This completes the proof of Lemma 4.\\

{\bf Acknowledgment} The author would like to express his sincere thanks to
the referee for his/her careful reading of the original
manuscript and giving of many kind suggestions to rewrite it. If this
revised manuscript has turned out to be readable, it is due to his/her advices.\\

\noindent {\bf \Large Appendix }\medskip

\noindent{\bf  A. The Sobolev imbedding theorem}\medskip

For the sake of self-containedness, we prove the Sobolev imbedding theorem
for our framework. (The statement is well-known if $N$
is an integer.) Let $y \in C^{\infty}[0,1]$ and $m \in \mathbb{N}$, we denote
$$\langle y\rangle_m^*:=\|(-\triangle)^my\|_{\mathfrak{X}}.$$
For $y\in \mathfrak{X}=L^2((0,1), x^{\frac{N}{2}-1}dx)$ we have the expansion $\displaystyle
y=\sum_{n=1}^{\infty}c_n\phi_n$, 
where
$(\phi_n)_n$ is the orthonormal system of eigenfunctions of
the operator $T=-\triangle$ with the Dirichlet boundary condition
at $x=1$. Then, for $m\in\mathbb{N}$ and for
$y \in C^{\infty}[0,1)$, we have
$$(-\triangle)^my(x)=\sum_{n=1}^{\infty}c_n\lambda_n^m\phi_n(x), \quad
\langle y\rangle_m^*=\Big(\sum_n
|c_n|^2\lambda_n^{2m}
\Big)^{1/2}. $$\\

\noindent{\bf Lemma A.1.}
Let $j_{\nu,n}$ be the $n$-th positive zero of the Bessel
function $J_{\nu}$, where $\nu=\frac{N}{2}-1$. Then we have
$$\lambda_n=(j_{\nu,n}/2)^2 \sim \frac{\pi^2}{4}n^2 \ \ \mbox{as }n\rightarrow\infty.$$

{\bf Proof}
By the Hankel's asymptotic expansion (see \cite{Watson}), the zeros of $J_{\nu}$ can be determined by
the relation
$$\tan\Big(r-\Big(\frac{\nu}{2}+\frac{1}{4}\Big)\pi\Big)=\frac{2}{\nu^2-\frac{1}{4}}r(1+O(r^{-2})).$$
Then we see
$$j_{\nu,n}=\Big(n_0+n+\frac{\nu}{2}+\frac{3}{4}\Big)\pi +
O\Big(\frac{1}{n}\Big) \ \ \mbox{as }n\rightarrow\infty,$$
for some $n_0\in\mathbb{Z}$. \hfill$\blacksquare$

\noindent{\bf Lemma A.2.}  There is a constant $C=C(N)$ such that
$$|\phi_n(x)|\leq C n^{\frac{N-1}{2}} \ \ \mbox{for }0\leq x\leq 1.$$

{\bf Proof}
Note that $\phi_n(x)$ is a normalization of $\Phi_{\nu}(\lambda_n x)$,
where
$$\Phi_{\nu}\Big(\frac{r^2}{4}\Big)=J_{\nu}(r)\Big(\frac{r}{2}\Big)^{-\nu}.$$
Since $|\Phi_{\nu}(x)|\leq C$ for $0\leq x<\infty$, it is sufficient to estimate
$\|\Phi_{\nu}(\lambda_nx)\|_{\mathfrak{X}}$. Using
the Hankel's asymptotic
expansion in the form
\begin{align*}
J_\nu(r)=&\sqrt{\frac{2}{\pi r}}\Big(\cos\Big(r-\frac{\nu}{2}\pi-\frac{\pi}{4}\Big)
\Big(1+O\Big(\frac{1}{r^2}\Big)\Big)+\\
-&\frac{1}{r}\sin\Big(r-\frac{\nu}{2}\pi-\frac{\pi}{4}\Big)
\Big(\frac{\nu^2-\frac{1}{4}}{2}+O\Big(\frac{1}{r^2}\Big)\Big)\Big),
\end{align*}
 we see that
\begin{align*}
\|\Phi_{\nu}(\lambda_nx)\|_{\mathfrak{X}}^2=&
(\lambda_n)^{-\nu-1}\int_0^{j_{\nu,n}}J_{\nu}(r)^2rdr
= (\lambda_n)^{-\nu-1}\Big(\frac{1}{\pi}j_{\nu,n}+O(1)\Big) \\
=&(\lambda_n)^{-\nu-1}\cdot \frac{2}{\pi}(\lambda_n^{1/2}+O(1))
\sim
\frac{2}{\pi}(\lambda_n)^{-\nu-\frac{1}{2}}.
\end{align*}
Then Lemma A.1 implies that
$$\|\Phi_{\nu}(\lambda_nx)\|_{\mathfrak{X}}^{-1}\sim \mbox{Const.}n^{\nu+\frac{1}{2}}.$$
\hfill$\blacksquare$

\noindent{\bf Lemma A.3.} If $ y \in C_0^{\infty}[0,1)$ and $0\leq j\leq m$, then
$\langle y\rangle_j^*\leq\langle y\rangle_m^*.$

{\bf Proof}
For
$y=\sum c_n\phi_n$, we have
\begin{align*}
(\langle y\rangle_j^*)^2=&
\sum |c_n|^2\lambda_n^{2j} =(\lambda_1)^{2j}\sum |c_n|^2(\lambda_n/\lambda_1)^{2j} \\
\leq& (\lambda_1)^{2j}\sum|c_n|^2(\lambda_n/\lambda_1)^{2m}=
\lambda_1^{2j-2m}(\langle y\rangle_m^*)^2.
\end{align*}
According to \cite[Section 15-6 (p.208)]{Watson}, we know that $j_{\nu,1}$ is an increasing function of
$\nu>0$ and $j_{\frac{1}{2}, 1}=\pi$. Therefore,
$\lambda_1\geq (\pi/2)^2 >1$ for $N\geq 2$ and which implies $\langle y\rangle_j^*\leq
\langle y\rangle_m^*$. \hfill$\blacksquare$

\noindent{\bf Lemma A.4.} If $2\sigma>N/2$, then
there is a constant $C
=C(\sigma,N)$ such that
$$\|y\|_{L^{\infty}}\leq C\langle y\rangle_{\sigma}^*$$
for any $y\in C_0^{\infty}[0,1)$.

{\bf Proof}
Let $y=\sum c_n\phi_n(x)$. Then Lemmas A.1 and A.2 imply
that
$$
|y(x)|\leq
\sum |c_n||\phi_n(x)|
\leq
C\sum |c_n|n^{\frac{N-1}{2}} 
\leq C\sqrt{\sum|c_n|^2\lambda_n^{2\sigma}}
\sqrt{\sum n^{N-4\sigma-1}}.
$$
Since $N-4\sigma<0$, the last term in the above inequality is finite. Therefore we get the required estimate. 
\hfill$\blacksquare$

Now, for $R>0$, we denote by $\mathfrak{X}(0,R)$ the Hilbert space of functions $y(x)$
of $0\leq x\leq R$ endowed with the inner product
$$(y_1|y_2)_{\mathfrak{X}(0,R)}=
\int_0^Ry_1(x)\overline{y_2(x)}x^{\frac{N}{2}-1}dx. $$
Moreover, for $m\in\mathbb{N}$, we denote by $\mathfrak{X}^{2m}(0,R)$ the space of functions $y(x)$ of
$0\leq x\leq R$ for which the derivatives
$(-\triangle)^jy \in \mathfrak{X}$ exist in the sense of distribution for $0\leq j\leq m$. And we use the norm
$$\|y\|_{\mathfrak{X}^{2m}(0,R)}:=\Big(
\sum_{0\leq j\leq m}\|(-\triangle)^jy\|_{\mathfrak{X}(0,R)}^2
\Big)^{1/2}. $$
Let us denote by $\mathfrak{X}_0^{2m}(0,R)$ the closure of $C_0^{\infty}[0,R)$ in the space
$\mathfrak{X}^{2m}(0,R)$. There is a continuous linear extension $\Psi: \mathfrak{X}^{2m}(0,1) \rightarrow
\mathfrak{X}_0^{2m}(0,2)$ such that
$$\|y\|_{\mathfrak{X}^{2m}(0,1)}\leq \|\Psi y\|_{\mathfrak{X}^{2m}(0,2)}
\leq C\|y\|_{\mathfrak{X}^{2m}(0,1)}.$$
See \cite[ p.186, Theorem 3.11]{Mizohata}, keeping in mind
Propositions 6,7. Then,
 by Lemmas A.3 and A.4, the Sobolev imbedding theorem holds for
$y \in \mathfrak{X}_0^{2\sigma}(0,2)$. That is, if $2\sigma >N/2$, there is a constant $C$
such that 
$\|y\|_{L^{\infty}}\leq C\|y\|_{\mathfrak{X}^{2\sigma}(0,2)} $
for $y \in \mathfrak{X}_0^{2\sigma}(0,2)$. Thus the same imbedding theorem holds for
$y \in C^{\infty}[0,1] \subset \mathfrak{X}^{2\sigma}(0,1)$
through the above extension. The conclusion is that, if $2\sigma >N/2$, there is a constant
$C=C(\sigma,N)$ such that
$\|y\|_{L^{\infty}}\leq C\sup_{0\leq j \leq \sigma}\|(-\triangle)^jy
\|_{\mathfrak{X}} $
for any $y \in C^{\infty}[0,1]$.\\

\noindent{\bf  B. Nirenberg-Moser type inequalities}\medskip

Let us prove Propositions 4, 5, 6 and 7.\\

\noindent{\bf Proof of Proposition 4}\medskip

First, it is easy to verify the formula
$$\dot{D}^kDy(x)=x^{-\frac{N+k}{2}}\int_0^x\dot{D}^k\triangle y(x')
(x')^{\frac{N+k}{2}-1}dx',\eqno(B.1)
$$
where $k\in \mathbb{N}$,
$$\dot{D}:=\sqrt{x}\frac{d}{dx}\quad\mbox{and} \quad D:=\frac{d}{dx}. $$
Since $\triangle =\dot{D}^2+\frac{N-1}{2}D$, (B.1) implies
$$|\dot{D}^kDy|_0\leq \frac{2}{N+k}|\dot{D}^{k+2}y|_0+\frac{N-1}{N+k}|\dot{D}^kDy|_0. $$
Here and hereafter $|\cdot|_0$ stands for $\|\cdot\|_{L^{\infty}}$. Thus we have
$$|\dot{D}^kDy|_0\leq \frac{2}{k+1}|\dot{D}^{k+2}y|_0.$$
Repeating this estimate, we get
$$|\dot{D}^kD^jy|_0\leq \Big(\frac{2}{k+1}\Big)^j|\dot{D}^{k+2j}y|_0. \eqno(B.2)$$
On the other hand, since
$\dot{D}^2=\triangle -\frac{N-1}{2}D$ and $D\triangle -\triangle D=D^2$, we have
$$\dot{D}^{2\mu}=\sum_{k=0}^{\mu}C_{k\mu}\triangle^{\mu-k}D^k \eqno(B.3)$$
with some constants $C_{k\mu}=C(k,\mu,N)$.
Then it follows from (B.3) and Proposition 3 that
$$|\dot{D}^{2\mu}D^jy|_0\leq C|\triangle^{\mu+j}y|_0. \eqno(B.4)$$
Since
$$\triangle =\dot{D}^2+\frac{N-1}{2}D\quad\mbox{and}\quad D\dot{D}^2-\dot{D}^2D=D^2, $$
it is easy to see that there are constants
$C_{km}=C(k,m,N)$ such that
$$\triangle^m=\sum_{k=0}^m
C_{km}\dot{D}^{2(m-k)}D^k. \eqno(B.5)$$
Applying the Leibnitz' rule to $D$ and $\dot{D}$, we see
$$\triangle^m(f\cdot g)=
\sum C_{k\ell jm}(\dot{D}^{2(m-k)-\ell}D^{k-j}f)\cdot(\dot{D}^{\ell}D^jg) \eqno(B.6)$$
with some constants $C_{k\ell jm}$. The summation is taken for
$0\leq j\leq k\leq m, 0\leq \ell \leq 2(m-k)$. By estimating
each term of the right-hand side of
(B.6), we can obtain the assertion of Proposition 4. In fact, we consider
the term
$$(\dot{D}^{\ell'}D^{j'}f)\cdot(\dot{D}^{\ell}D^jg)$$
provided that
$\ell'+\ell +2(j'+j)=2m$. By (B.2) and (B.4) we have
\begin{align*}
|\dot{D}^{\ell}D^jg|_0\leq
C|\dot{D}^{\ell +2j}g|_0
\leq  C'|\dot{D}^{2m}g|_0^{\frac{\ell +2j}{2m}}|g|_0^{1-\frac{\ell+2j}{2m}}
\leq  C''|\triangle^mg|_0^{\frac{\ell+2j}{2m}}|g|_0^{1-\frac{\ell+2j}{2m}}
\end{align*}
for some positive constants $C$, $C'$ and $C''$.
Here we have used the Nirenberg interpolation for $\dot{D}=\partial/\partial \xi$,
where $x=\xi^2/4$. The same
estimate holds for
$|\dot{D}^{\ell'}D^{j'}f|_0$. Therefore we have
\begin{align*}
|(\dot{D}^{\ell'}D^{j'}f)\cdot(\dot{D}^{\ell}D^jg)|_0 \leq &
 C|\triangle^mf|_0^{\frac{\ell'+2j'}{2m}}|f|_0^{1-\frac{\ell'+2j'}{2m}}|\triangle^mg|_0^{\frac{\ell+2j}{2m}}
|g|_0^{1-\frac{\ell+2j}{2m}} \\
\leq & C(|\triangle^mf|_0|g|_0+|f|_0|\triangle^mg|_0),
\end{align*}
since $X^{\theta}Y^{1-\theta}\leq X+Y$. \\

\noindent{\bf Proof of Proposition 5}\medskip

Suppose $F(x,y)$ is a smooth function of $x$ and $y$. Let us consider
the composed function $U(x):=F(x,y(x))$. We claim that
$$|\triangle^mU|_0\leq C(1+|y|_m^*) $$
provided that $|y|_0\leq M$.
 In fact,
$$\triangle^mU=\sum C_{km}\dot{D}^{2(m-k)}D^k U $$
consists of several terms of the following form:
$$\Big(\dot{D}_x^K\Big(\frac{\partial}{\partial y}\Big)^LD_x^k\Big(\frac{\partial}{\partial y}\Big)^{\ell}F\Big)\cdot
(\dot{D}^{K_1}y)\cdots(\dot{D}^{K_L}y)\cdot(\dot{D}^{\mu_1}D^{k_1}y)\cdots
(\dot{D}^{\mu_{\ell}}D^{k_{\ell}}y), $$
where
$$k+k_1+\cdots+k_{\ell}=\kappa,$$ $$ K+K_1+\cdots K_L+\mu_1+\cdots+\mu_{\ell}=2(m-\kappa).$$
Therefore
$$K_1+\cdots+K_L+(\mu_1+2k_1)+\cdots(\mu_{\ell}
+2k_{\ell}) \leq 2m. $$
Applying the Nirenberg interpolation to $\dot{D}$ and using (B.4),
we have
$$|\dot{D}^{K_1}y|_0\leq C(|y|_m^*)^{\frac{K_1}{2m}}|y|_0^{1-\frac{K_1}{2m}}. $$
Similarly,
$$|\dot{D}^{\mu_1}
D^{k_1}y|_0\leq C|\dot{D}^{\mu_1+2k_1}y|_0\leq C'(|y|_m^*)^{\frac{\mu_1+2k_1}{2m}}
|y|_0^{1-\frac{\mu_1+2k_1}{2m}}, $$
and so on. Then our claim follows obviously. \\

\textbullet\  We note that by (B.2), (B.4) and (B.5) we have
$$\frac{1}{C}|\dot{D}^{2j}f|_0\leq |\triangle^jf|_0
\leq C|\dot{D}^{2j}f|_0.
\eqno(B.7)$$ \\

\noindent{\bf Proof of Proposition 6}\medskip

It can be verified that
$$\triangle^m(\alpha y)=
\sum_{j+k=m}(\alpha_{1k}^{(m)}\check{D}\triangle^jy+
\alpha_{0k}^{(m)}\triangle^jy), $$
where
$\alpha_{1k}^{(m)}$ and $\alpha_{0k}^{(m)}$ are determined by the recurrence formula
\begin{align*}
\alpha_{1k}^{(m+1)}&=\alpha_{1k}^{(m)}+
(\triangle -(N-2)D)\alpha_{1,k-1}^{(m)}+2D\alpha_{0,k-1}^{(m)},\\
\alpha_{0k}^{(m+1)}&=(1+2\check{D})\alpha_{1k}^{(m)}+
\alpha_{0k}^{(m)}+\triangle\alpha_{0,k-1}^{(m)},
\end{align*}
starting from
$\alpha_{10}^{(0)}=0, \alpha_{00}^{(0)}=\alpha$. 
Here we have used the convention $\alpha_{1k}^{(m)}=\alpha_{0k}^{(m)}
=0$ for $k<0$ or $k>m$. Of course $\alpha_{10}^{(m)}=0$ for any $m$.
Therefore we see that $\|\triangle^m(\alpha y)\|_0\leq C \|y\|_{2m}$.

Differentiating the formula, we get
$$
\dot{D}\triangle^m(\alpha y) =
\sum_{j+k=m}(\dot{\alpha}_{2k}^{(m)}\triangle^{j+1}y+\dot{\alpha}_{1k}^{(m)}\dot{D}\triangle^jy +
\dot{\alpha}_{0k}^{(m)}\triangle^jy),
$$
where
$$
\dot{\alpha}_{2k}^{(m)}=\sqrt{x}\alpha_{1k}^{(m)}, \quad
\dot{\alpha}_{1k}^{(m)}=\Big(-\frac{N}{2}+1+\check{D}\Big)\alpha_{1k}^{(m)}+
\alpha_{0k}^{(m)}, \quad
\alpha_{0k}^{(m)}=\dot{D}\alpha_{0k}^{(m)}.
$$
It is clear that $\|\dot{D}\triangle^m(\alpha y)\|_0
\leq C\|y\|_{2m+1}$, since $\dot{\alpha}_{20}^{(m)}=0$
for any $m$. \hfill$\blacksquare$\\

\noindent{\bf Proof of Proposition 7}\medskip

It is sufficient to prove that
$$\|\dot{D}y\|\leq C(\|y\|+\|\triangle y\|), $$
where and hereafter we denote $\|\cdot\|=\|\cdot\|_{\mathfrak{X}}$.

If $w$ satisfies the Dirichlet boundary condition $w(1)=0$, then
$$\|\dot{D}w\|^2=(-\triangle w\ |\ w)\leq \|\triangle w\|\|w\|.$$
Therefore
we have
$$\|\dot{D}y\|^2\leq \|\triangle y\|(\|y\|+|y(1)|).$$
On the other hand we have
$$\sqrt{\frac{2}{N}}|y(1)|\leq \|y\|+\sqrt{\frac{2}{N-2}}\|\dot{D}y\|.$$
In fact, since
$$y(1)=y(x)+\int_x^1\frac{1}{\sqrt{x}}\dot{D}y(x')dx',$$
we have
$$|y(1)|^2\leq |y(x)|^2 +\frac{2}{N-2}\|\dot{D}y\|^2x^{-\frac{N}{2}+1}$$
for $x>0$. Integrating this, we get the above estimate of $|y(1)|$.
Hence we have, for any $\epsilon>0$,
\begin{align*}
\|\dot{D}y\|^2&\leq C\|\triangle y\|(\|y\|+\|\dot{D}y\|) \\
&\leq C\Big(\frac{1}{2\epsilon}\|\triangle y\|^2+\frac{\epsilon}{2}(\|y\|+\|\dot{D}y\|)^2\Big) \\
&\leq C\big(\frac{1}{2\epsilon}\|\triangle y\|^2+\epsilon \|y\|^2+
\epsilon \|\dot{D}y\|^2\Big).
\end{align*}
Taking $\epsilon$ to be small, we get the desired estimate. \hfill$\blacksquare$\\

\noindent{\bf  C. Existence of the smooth solution to the linear wave equation}\medskip

Let us give a proof of the existence of the smooth solution
to the initial boundary value problem (IBP):
$$
\frac{\partial^2h}{\partial t^2}+\mathcal{A}h=g(t,x), \quad h|_{x=1}=0, \quad
h|_{t=0}=\frac{\partial h}{\partial t}\Big|_{t=0}=0.
$$
We assume that $g(t,x)=0$ for $0\leq t\leq \tau_1$.\\

{\bf Existence.} The existence of the solution can be proved by applying the Kato's theory
developed in \cite{Kato1970}. In fact, we consider the closed operator
$$\mathfrak{A}(t)=\begin{bmatrix} 0 & -1\\
\mathcal{A}(t) & 0 \end{bmatrix}
$$
in $\mathfrak{H}:=\mathfrak{X}_0^1\times\mathfrak{X}$ densely defined on
$$\mathcal{D}(\mathfrak{A}(t))=\mathfrak{G}:=\mathfrak{X}_{(0)}^2\times\mathfrak{X}^1.$$
Here $\mathfrak{X}=L^2((0,1);x^{\frac{N}{2}-1}dx), \mathfrak{X}^1=\{y\in\mathfrak{X} |
\dot{D}y\in\mathfrak{X}\}, \mathfrak{X}_0^1=\{y\in\mathfrak{X}^1 | y|_{x=1}=0\}, 
\mathfrak{X}^2=\{y\in\mathfrak{X}^1 | \triangle y\in \mathfrak{X}\}$ and
$\mathfrak{X}_{(0)}^2=\mathfrak{X}^2\cap\mathfrak{X}_0^1=\{y\in\mathfrak{X}^2 |
y|_{x=1}=0\}$. The problem (IBP) is equivalent to
$$\frac{du}{dt}+\mathfrak{A}(t)u=\mathfrak{f}(t), \qquad u|_{t=0}=0,$$
with
$$ u=\begin{bmatrix} h \\
\frac{\partial h}{\partial t}\end{bmatrix} \quad \mbox{and}\quad \mathfrak{f}(t)=\begin{bmatrix} 0 \\
g(t,\cdot)\end{bmatrix}.$$
We can write
$$\mathcal{A}(t)y=-x^{-\frac{N}{2}+1}\frac{d}{dx}ax^{\frac{N}{2}}\frac{dy}{dx}+b\check{D}y
+cy,$$
where
$$a=b_2,\qquad b=b_1+Db_2,\qquad c=b_0.$$
Then
$$(\mathcal{A}(t)y|v)_{\mathfrak{X}}=(a(t)\dot{D}y|\dot{D}v)_{\mathfrak{X}}+((b\check{D}+c)y|v)_{\mathfrak{X}}$$
for $y \in \mathfrak{X}_{(0)}^2
$ and $v\in\mathfrak{X}_0^1$. The inner product
$$(y|v)_t=(a(t)\dot{D}y|\dot{D}v)_{\mathfrak{X}}+(y|v)_{\mathfrak{X}}$$
introduces an equivalent norm $\|\cdot\|_t$ in $\mathfrak{X}_0^1$ provided
that $|1-a|\leq 1/2$, $\|a\|_{L^{\infty}}, \|b\|_{L^{\infty}}, \|c\|_{L^{\infty}}
\leq M_0.$ Then we have
$$-(\mathfrak{A}(t)u|u)_{\mathfrak{H}_t}=(u_2|u_1)_{\mathfrak{X}}-
((b\check{D}+c)u_1|u_2)_{\mathfrak{X}} \leq \beta \|u||_{\mathfrak{H}_t}^2,$$
where
\begin{align*}
(u|\phi)_{\mathfrak{H}_t}&=(u_1|\phi_1)_t+(u_2|\phi_2)_{\mathfrak{X}}\\
&=(a(t)\dot{D}u_1|\dot{D}\phi_1)_{\mathfrak{X}}+
(u_1|\phi_1)_{\mathfrak{X}}+(u_2|\phi_2)_{\mathfrak{X}}
\end{align*}
and $\beta$ depends only upon $M_0$. $\|u\|_{\mathfrak{H}_t}=\sqrt{(u|u)_{\mathfrak{H}_t}}$ is
equivalent to $\|u\|_{\mathfrak{H}}$ and depends on $t$ smoothly in the sense of
\cite[Proposition 3.4]{Kato1970}. From the above estimate it follows that
$\mathfrak{A}(t)$ is a quasi-accretive generator in the norm $\|\cdot\|_{\mathfrak{H}_t}$.

In fact the following argument is standard: the equation
$$(\lambda+\mathfrak{A}(t))u=f$$
is reduced to an elliptic equation
$$(\lambda^2+\mathcal{A}(t))u_1=\lambda f_1+f_2,$$
which admits a solution $u_1\in \mathfrak{X}_{(0)}^2$ for given
$f_3:=\lambda f_1+f_2\in\mathfrak{X}$, provided that
$\lambda^2>\|b\|_{L^{\infty}}^2+\|c\|_{L^{\infty}}+\frac{1}{4}$; then 
$$Q[u]:=\lambda^2\|u\|_{\mathfrak{X}}^2+
(a(t)\dot{D}u|\dot{D}u)_{\mathfrak{X}}+((b\check{D}+c)u|u)_{\mathfrak{X}}\geq \frac{1}{4}\|u\|_{\mathfrak{X}^1}^2,$$
and for given $f_3\in\mathfrak{X}$ there is a $u_1\in \mathfrak{X}_0^1$ such that
$Q(u_1,v)=(f|v)_{\mathfrak{X}}$ for any $v\in \mathfrak{X}_0^1$; thus
$(\lambda+\mathfrak{A}(t))^{-1}\in\mathcal{B}(\mathfrak{H})$ and $\|(\lambda+\mathfrak{A}(t))^{-1}\|_{\mathcal{B}(\mathfrak{H}_t)}
\leq(\lambda-\beta)^{-1}$.

Therefore by \cite[Proposition 3.4]{Kato1970}, $(\mathfrak{A}(t))_t$ is a stable
family of generators. Hence by \cite[Theorem 7.1 and 7.2]{Kato1970}, we can claim that
there exists a solution $u\in C^1([0,T]; \mathfrak{H})\cap C([0,T];\mathfrak{G})$, which
gives
the desired solution
$h\in C^2([0,T];\mathfrak{X})\cap C^1([0,T];\mathfrak{X}^1)\cap
C([0,T];\mathfrak{X}^2)$, since $g\in C^{\infty}([0,T]\times[0,1])$.

{\bf Regularity.} We want to show that $h \in C^{\infty}([0,T]\times[0,1])$. 
To do so, we apply the Kato's theory developed in \cite[Section 2]{Kato1976}. 
We consider the spaces
\begin{align*}
&\hat{\mathfrak{H}}=\hat{\mathfrak{H}}_0=\mathfrak{X}_0^1\times\mathfrak{X}\times\mathbb{R}, \\
&\hat{\mathfrak{H}}_j=\mathfrak{X}_{(0)}^{j+1}\times\mathfrak{X}^j\times\mathbb{R}, \\
&\hat{\mathfrak{G}}=\hat{\mathfrak{G}}_1=\mathfrak{X}_{(0)}^2\times \mathfrak{X}_0^1\times\mathbb{R}, \\
&\hat{\mathfrak{G}}_j=\hat{\mathfrak{G}}\cap\hat{\mathfrak{H}}_j=\mathfrak{X}_{(0)}^{j+1}\times
\mathfrak{X}_{(0)}^j\times\mathbb{R}.
\end{align*}
Here $\mathfrak{X}^k=\{y | \|y\|_k=(\sum_{0\leq\ell\leq k}\langle y\rangle_{\ell}^2)^{1/2}<\infty\},$
and $\mathfrak{X}_{(0)}^k=\mathfrak{X}^k\cap \mathfrak{X}_0^1$.
Introducing the closed operator
$$
\hat{\mathfrak{A}}(t)=\begin{bmatrix} 0 & -1 & 0 \\
\mathcal{A}(t) & 0 & -g(t) \\
0 & 0 & 0
\end{bmatrix}$$
in $\hat{\mathfrak{H}}$
densely defined on
$$\mathcal{D}(\hat{\mathfrak{A}}(t))=\hat{\mathfrak{G}},$$
we can convert (IBP) to
$$\frac{du}{dt}+\hat{\mathfrak{A}}(t)u=0,\qquad u|_{t=0}=\phi_0, $$
where
$$\phi_0=\begin{bmatrix}0\\
0\\
1\end{bmatrix}.$$
Since $g \in C^{\infty}$, the stability of $(\hat{\mathfrak{A}}(t))_t$ is
reduced to that of $(\mathfrak{A}(t))_t$ by the perturbation theorem (\cite[Proposition 1.2]{Kato1976}).
Therefore $(\hat{\mathfrak{A}}(t))_t$ is a stable family of generators in
$\hat{\mathfrak{H}}$. Since the coefficients of the
differential operator $\mathcal{A}$ are in $C^{\infty}$, we see
$\mathcal{D}(\hat{\mathfrak{A}}(t))\cap\hat{\mathfrak{H}}_1=
\hat{\mathfrak{G}}$ and
$$\frac{d^k}{dt^k}\hat{\mathfrak{A}}(t)\in L^{\infty}([0,T];
\mathcal{B}(\hat{\mathfrak{G}}_{j+1}, \hat{\mathfrak{H}}_j))$$
for all $j,k$. Moreover we have `ellipticity', i.e., 
for each $t,j$ $u\in \mathcal{D}(\hat{\mathfrak{A}}(t))$ and
$\hat{\mathfrak{A}}(t)u\in \hat{\mathfrak{H}}_j$ implies
$u\in \hat{\mathfrak{H}}_{j+1}$ with
$$\|u\|_{\hat{\mathfrak{H}}_{j+1}}\leq C(\|\hat{\mathfrak{A}}(t)u\|_{\hat{\mathfrak{H}}_j}+
\|u\|_{\hat{\mathfrak{H}}}).
$$
In fact this condition is reduced to the fact that if
$y \in \mathfrak{X}^2$ and $\mathcal{A}(t)y \in \mathfrak{X}^j$ then
$y \in \mathfrak{X}^{j+2}$ and
$$\|y\|_{j+2}\leq C(\|\mathcal{A}(t)y\|_j+\|y\|_1).$$
See Proposition 8. Thus we can apply
\cite[Theorem 2.13]{Kato1976}, say, if
$\phi_0\in D_m(0)$, then the solution $u$ satisfies
$$u \in \bigcap_{j+k=m}C^k([0,T]; \hat{\mathfrak{H}}_j), $$
which implies
$$h \in \bigcap_{j+k=m}C^k([0,T];\mathfrak{X}_{(0)}^{j+1}).$$
Recall $\phi_0=(0,0,1)^T$ and the space of compatibility
$D_m(0)$ is characterized by
\begin{align*}
D_0(0)&=\hat{\mathfrak{H}},  \quad
S^0(0)=I, \\
D_{j+1}(0)&=\{\phi\in D_j(0) |
S^k(0)\phi \in \hat{\mathfrak{G}}_{j+1-k}, 0\leq k\leq j\}, \\
S^{j+1}(0)\phi&=-
\sum_{k=0}^j
\binom{j}{k}\Big(\frac{d}{dt}\Big)^{j-k}\hat{\mathfrak{A}}(0)S^k(0)\phi.
\end{align*}
See \cite[(2.40), (2.41)]{Kato1976}. Since $g=0$ for $0\leq t\leq \tau_1$,
we have
$$\Big(\frac{d}{dt}\Big)^{n}\hat{\mathfrak{A}}(0)=
\begin{bmatrix}0&0&0\\
\Big(\frac{d}{dt}\Big)^n\mathcal{A}(0) & 0 & 0\\
0&0&0
\end{bmatrix}.$$
Thus
it is easy to see $\phi_0\in D_j(0)$ and $S^j(0)\phi_0=0$ for 
$j\geq 1$ inductively on $j$. (Note that $S^0(0)\phi_0=\phi_0$.)
Hence for any positive integer $m$ we have $\phi_0\in D_m(0)$ and 
obtain the desired regularity of the solution $h$.\\

\noindent{\bf D. Proof of Proposition 8}\\

By tedious calculations we have
$$
[\triangle^m, \mathcal{A}]y :=\triangle^m\mathcal{A}y-\mathcal{A}\triangle^my 
=\sum_{j+k=m}(b_{1k}^{(m)}\check{D}\triangle^jy+b_{0k}^{(m)}\triangle^jy),
$$
where $\check{D}=xd/dx$ and
\begin{align*}
b_{10}^{(m)}&=-2mDb_2, \\
b_{00}^{(m)}&=-m((2m-1)\triangle +(m-1)(1-N)D)b_2+
m(1+2\check{D})b_1,
\end{align*}
where $D=d/dx$ and $b_{1k}^{(m)}, b_{0k}^{(m)}, k\geq 1$ are determined by 
$$ b_{11}^{(1)}=2Db_0+(\triangle -(N-2)D)b_1,\quad
 b_{01}^{(1)}=\triangle b_0, 
$$
and the recurrence formula
\begin{align*}
b_{1k}^{(m+1)}&=b_{1k}^{(m)}+(\triangle-(N-2)D)b_{1,k-1}^{(m)}
+2Db_{0,k-1}^{(m)}\quad\mbox{for}\quad k\geq 2, \\
b_{11}^{(m+1)}&=b_{11}^{(m)}-4m^2(\triangle+\frac{3-N}{2}D)Db_2 + \\
&+((4m+1)\triangle -(2mN-6m+N-2)D)b_1+2Db_0, \\
b_{0k}^{(m+1)}&= b_{0k}^{(m)}+(1+2\check{D})b_{1k}^{(m)}+\triangle b_{0,k-1}^{(m)}
\quad\mbox{for}\quad k\geq 2, \\
b_{01}^{(m+1)}&= b_{01}^{(m)}-m\triangle((2m-1)\triangle + (m-1)(1-N)D)b_2+\\
&+m(3+2\check{D})\triangle b_1+\triangle b_0+(1+2\check{D})b_{11}^{(m)}.
\end{align*}

We have used the following calculus formula:

\begin{align*}
D\check{D}&=\triangle -\Big(\frac{N}{2}-1\Big)D, \qquad
\triangle\check{D}-\check{D}\triangle =\triangle, \\
\triangle (Q\check{D}P)&=Q\check{D}\triangle P+(1+2\check{D})Q\cdot\triangle P+
(\triangle-(N-2)D)Q\cdot\check{D}P, \\
\triangle(QP)&=Q\triangle P+2(DQ)\check{D}P+(\triangle Q)P.
\end{align*}

Then it follows that
$$\|b_{0k}^{(m)}\|_{L^{\infty}}\leq C|\vec{b}|_{2k+3}
\quad \|b_{1k}^{(m)}\|_{L^{\infty}}\leq C|\vec{b}|_{2k+2} $$
and therefore
$$\|[\triangle^m, \mathcal{A}]y\|\leq
C\sum_{j+k=m}(|\vec{b}|_{2k+2}\|y\|_{2j+1}
+|\vec{b}|_{2k+3}\|y\|_{2j}). $$

Since
$\triangle^m[\triangle,\mathcal{A}]=[\triangle^{m+1}, \mathcal{A}]-
[\triangle^m, \mathcal{A}]\triangle, $
it follows that
$\|\triangle^m[\triangle,\mathcal{A}]y\|\leq C A_m,$
where
$$A_m:=\sum_{j+k=m+1}
(|\vec{b}|_{2k+2}\|y\|_{2j+1}+|\vec{b}|_{2k+3}\|y\|_{2j}). $$

{\bf Remark}\  This estimate is very rough and may be far from the best possible.
But it is enough for our purpose. To derive this estimate, we have used the following
observations:

Let $\mathfrak{M}_k$ denote the set of all functions of the form
$$\sum_{\alpha=2,1,0}\sum_{i+j\leq k}C_{\alpha ij}\triangle^iD^jb_{\alpha},$$
$C_{\alpha ij}$ being constants. Then it can be shown that $\triangle f$, $Df$
and $D\check{D}f(=\triangle f+(-\frac{N}{2}+1)Df)$
 belong to $\mathfrak{M}_{k+1}$ if $f$ belongs to $\mathfrak{M}_k$.
Using this, we can claim inductively that $b_{1k}^{(m)}\in
\mathfrak{M}_{k+1}$ and $b_{0k}^{(m)}\in\mathfrak{M}_{k+1}+\check{D}\mathfrak{M}_{k+1}$
for any $m$, $k\leq m+1$. Note that
$\|f\|_{L^{\infty}}\leq C|\vec{b}|_{2k}$
and $\|\check{D}f\|_{L^{\infty}}\leq\|\dot{D}f\|_{L^{\infty}}\leq C|\vec{b}|_{2k+1}$ if $f\in\mathfrak{M}_k$.
(See Proposition 3 and Appendix B, (B.7).
Also see (B.3), keeping in mind that
$\triangle=\dot{D}^2+\frac{N-1}{2}D$.)\hfill$\square$

Differentiating $[\triangle^m, \mathcal{A}]y$, we get
$$\dot{D}[\triangle^m, \mathcal{A}]y=
\sum_{k+j=m}(\dot{b}_{2k}^{(m)}\triangle^{j+1}y+
\dot{b}_{1k}^{(m)}\dot{D}\triangle^jy+
\dot{b}_{0k}^{(m)}\triangle^jy), $$
where
$$
\dot{b}_{2k}^{(m)}=\sqrt{x}b_{1k}^{(m)}, \quad
\dot{b}_{1k}^{(m)}=(-\frac{N}{2}+1+\check{D})b_{1k}^{(m)}+b_{0k}^{(m)}, \quad
\dot{b}_{0k}^{(m)}=\dot{D}b_{0k}^{(m)}.
$$
Using
$$\dot{D}\triangle^m[\triangle, \mathcal{A}]=\dot{D}[\triangle^{m+1}, \mathcal{A}]-
\dot{D}[\triangle^m, \mathcal{A}]\triangle, $$
we have
$\|\dot{D}\triangle^m[\triangle, \mathcal{A}]y\|\leq
CA_m^{\sharp},$
where
$$A_m^{\sharp}:=\sum_{j+k=m+1}(|\vec{b}|_{2k+2}\|y\|_{2j+2}
+|\vec{b}|_{2k+3}\|y\|_{2j+1}+|\vec{b}|_{2k+4}\|y\|_{2j}).$$

Since $A_{m-1}\leq A_{m-1}^{\sharp}\leq 2A_m\leq 2A_{m}^{\sharp}$,
we can claim that
$$
\|[\triangle, \mathcal{A}]y\|_{2m} \leq CA_m, \quad
\|[\triangle, \mathcal{A}]y\|_{2m+1}\leq CA_m^{\sharp}.
$$

Now
$\triangle y=-\frac{1}{b_2}(\mathcal{A}y-b_1\check{D}y-b_0y) $
implies 
$\|\triangle y\|\leq C(\|\mathcal{A}y\|+\|y\|_1), $
and
$\|y\|_2\leq C(\|\mathcal{A}y\|+\|y\|_1). $

Moreover
\begin{align*}
\dot{D}\triangle y &=-\frac{1}{b_2}\Big(\dot{D}\mathcal{A}y+(-\dot{D}b_2+\sqrt{x}b_1)\triangle y + \\
&+ (-\frac{N}{2}+1+\check{D}b_1+b_0)\dot{D}y+(\dot{D}b_0)y\Big)
\end{align*}
implies
$\|\dot{D}\triangle y\|\leq C(\|\dot{D}\mathcal{A}y\|+\|y\|_2),$
and
$\|y\|_3\leq C(\|\mathcal{A}y\|_1+\|y\|_1). $

Using the estimates of $[\triangle,\mathcal{A}]$, we can show inductively that,
for $n\geq 2$,
$$\|y\|_{n+2}\leq C(\|\mathcal{A}y\|_n+\|y\|_1+K(n)), $$
where
\begin{align*}
K(n)=\begin{cases}
A_m &\mbox{for $n=2m+2$,} \\
A_m^{\sharp} &\mbox{for $n=2m+3$}
\end{cases}
\end{align*} 

By interpolation we have
$$
K(n)\leq
C(|\vec{b}|_2\|y\|_{n+1}+|\vec{b}|_{n+3}\|y\|). 
$$
This completes the proof of Proposition 8.

\par\vspace{30mm}

nuna adreso:

Tetu Makino

Department of Applied Mathematics, Faculty of Engineering,

Yamaguchi University,

Ube 755-8611, Japan

E-mail: makino@yamaguchi-u.ac.jp


\begin{thebibliography}{9}
\bibitem{Adams} R. A. Adams, Sobolev Spaces, Academic Press, 1975.
\bibitem{Coddington} E. A. Coddington and N. Levinson, 
Theory of Ordinary Differential Equations, McGraw-Hill, 1955
\bibitem{CS} D. Coutand and S. Shkoller, Well-posedness in smooth function spaces
for moving-boundary 1-d compressible
Euler equations in physical vacuum,
Comm. Pure Appl. Math., LXIV(2011), 328-366
\bibitem{DCL} Ding Xiaxi, Chen Guiqiang and Luo Peizhu,
Convergence of the fractional step Lax-Friedrichs scheme and
Godunov scheme for the isentropic system of
gas dynamics, Comm. Math. Phys., 121(1989), 63-84
\bibitem{Hamilton} R. Hamilton, The inverse function theorem of Nash and Moser, Bull. American Math. Soc.,
7 (1982), 65-222
\bibitem{Ikawa} M. Ikawa, Hyperbolic Partial Differential
Equations and Wave Phenomena (Translations of Math. Monographs, Vol.189),
AMS, Providence, Rhode Island, 2000.
\bibitem{JM} J.-H. Jang and N. Masmoudi, 
Well-posedness for compressible Euler equations with
physical vacuum singularity, Comm. Pure Appl.
Math., LXII(2009), 1327-1385
\bibitem{Kato1970} T. Kato, Linear evolution equations of ``hyperbolic" type,
J. Fac. Sci. Univ. Tokyo, Section I, 17(1970), 241-258.
\bibitem{Kato1976} T. Kato, Linear and quasi-linear equations of evolution
of hyperbolic type, in Hyperbolicity, CIME, II Ciclo, 1976, 125-191; 
reprinted by Springer, 2011
\bibitem{L} T.-P. Liu, Compressible flow with damping and vacuum,
Japan J. Appl. Math., 13(1996), 25-32
\bibitem{LY} T.-P. Liu and T. Yang, Compressible flow with vacuum
and physical singularity, Methods Appl. Anal., 31(2000), 223-237
\bibitem{M1986} T. Makino, On a local existence theorem for the evolution
of gaseous stars, in: Patterns and Waves,
ed. by T. Nishida, M. Mimura and H. Fujii, North-Holland,
1986, 459-479
\bibitem{M1989} T. Makino, Les solutions \`{a} support compact de l'\'{e}quation
du mouvement des atmosph\`{e}res d'\'{e}toiles, Japan J. Appl. Math.,
6(1989), 479-489
\bibitem{MT} T. Makino and S. Takeno, Initial boundary value problem for the spherically symmetric
motion of isentropic gas, Japan J. Indust. Appl. Math., 11(1994), 171-183
\bibitem{MUK} T. Makino, S. Ukai et S. Kawashima, 
Sur la solution \`{a} support compact de l'\'{e}quation d'Euler compressible,
Japan J. Appl. Math., 3(1986), 246-257
\bibitem{Mizohata} S. Mizohata, The Theory of Partial Differential
Equations, Cambridge University Press, 1973. 
\bibitem{Reed} M. Reed and B. Simon, Methods of Modern Mathematical
Physics, Vol.II: Fourier Analysis, Self-Adjointness, Academic Press, 1975
\bibitem{Sneddon} I. N. Sneddon, Fourier Transformations, NY, McGraw-Hill, 1951; NY, Dover,
1995.
\bibitem{Watson} G. N. Watson, A Treatise on the Theory of
Bessel Functions, Cambridge University Press, 1958.
\bibitem{Y} T. Yang, Singular behavior of vacuum states for compressible fluids,
Comput. Appl. Math., 190(2006), 211-231


\end{thebibliography}
\end{document}